\documentclass[leqno,12pt]{amsart}

\usepackage{amssymb}
\usepackage{amscd}              % for commutative diagrams
\usepackage[all]{xy}            % for use of xypic
\xyoption{curve}
\CompileMatrices

\newcommand{\Aut}{\text{Aut}}
\newcommand{\Bc}{{\mathcal{B}}}

\newcommand{\Bf}{{\mathfrak{B}}}

\newcommand{\C}{\mathbb{C}}

\newcommand{\Cc}{{\mathcal{C}}}
\newcommand{\Cf}{{\mathfrak{C}}}

\newcommand{\comp}{\circ}

\newcommand{\Ct}{\tilde{C_0}}                  
\newcommand{\divi}{\text{div}\, }                % divisor
\newcommand{\Dc}{{\mathcal{D}}}
\newcommand{\Dk}{\Delta^\kappa}

\newcommand{\eins}{{\mathbb{I}}}

\newcommand{\E}{{\mathcal{E}}}

\newcommand{\Et}{\text{$\tilde{\E}$}}

\newcommand{\F}{{\mathcal{F}}}
\newcommand{\Fl}{\text{\bf Fl}}
\newcommand{\G}{{\mathcal{G}}}

\newcommand{\GL}{\text{GL}}

\newcommand{\Gm}{\mathbb{G}_m}

\newcommand{\GVB}{\text{GVB}}
\newcommand{\GVBD}{\text{GVBD}}

\newcommand{\Hc}{{\mathcal{H}}}

\newcommand{\id}{\text{id}}

\newcommand{\injto}{\hookrightarrow}
\newcommand{\isomorph}{\cong}           
\newcommand{\isomto}{\overset{\sim}{\rightarrow}}  

\newcommand{\Id}{\text{Id}}

\newcommand{\Isom}{\text{Isom}}
\newcommand{\Isomto}{\overset{\sim}{\longrightarrow}}
\newcommand{\J}{{\mathcal{J}}}

\newcommand{\KGL}{\text{KGL}}
\newcommand{\KGLn}{\text{KGL}_n}

\newcommand{\Liset}{\text{Lis-\'et}}
\newcommand{\Ll}{{\mathcal{L}}}                  % line bundle L
\newcommand{\m}{{\mathfrak{m}}}                  % maximal ideal m
\newcommand{\M}{{\mathcal{M}}}

\newcommand{\ob}{\text{ob}}
\newcommand{\Ob}{\overline{\text{\bf O}}}
\newcommand{\Oo}{{\mathcal{O}}}                  % structure sheaf
\newcommand{\Oplus}{\bigoplus}
\newcommand{\PB}{\text{PB}}

\newcommand{\Pic}{\text{Pic}}
\newcommand{\pit}{\tilde{\pi}}

\newcommand{\Pp}{{\mathbb{P}}}                   % projective space P
\newcommand{\pr}{\text{pr}}
\newcommand{\Proj}{\text{Proj}\,}
\newcommand{\Qb}{\overline{Q}}

\newcommand{\res}{\text{res}}
   
\newcommand{\SA}{\text{SA}}
\newcommand{\Spec}{\text{Spec}\, }               % spectrum of a ring

\newcommand{\sign}{\text{sign}}
\newcommand{\SL}{\text{SL}}
\newcommand{\sln}{\mathfrak{sl}_n}
\newcommand{\SPB}{\text{SPB}}
\newcommand{\SVB}{\text{SVB}}

\newcommand{\tensor}{\otimes}

\newcommand{\Tensor}{\bigotimes}

\newcommand{\Thetat}{\tilde{\Theta}}
\newcommand{\Tk}{\Theta^\kappa}
\newcommand{\Tr}{\text{Tr}}

\newcommand{\To}{\longrightarrow}

\newcommand{\VB}{\text{VB}}

\newcommand{\X}{{\mathcal{X}}}

\newcommand{\Z}{{\mathbb{Z}}}                    % integers

\newtheorem{theorem}{Theorem}[section]
\newtheorem{proposition}[theorem]{Proposition}
\newtheorem{lemma}[theorem]{Lemma}

\theoremstyle{definition}
\newtheorem{definition}[theorem]{Definition}
\newtheorem{remark}[theorem]{Remark}

\newtheorem{sketch of proof}[theorem]{Sketch of proof}

\makeatletter
\def\overunderbraces #1#2#3{{%
 \baselineskip\z@skip \lineskip4\p@ \lineskiplimit4\p@
 \displaystyle  %% generate error if not in math mode!
 \setbox\z@\vbox{\ialign{&\hfil${}##{}$\hfil\cr
   \global\let\br\br@label #1\cr %upper labels
   \global\let\br\br@down #1\cr   %upper braces
   #2\cr % main line of the formula
 }}% finished partial alignment and \vbox
 \dimen@-\ht\z@ %   Measure height of partial alignment ..
 \setbox\z@\vbox{\ialign{&\hfil${}##{}$\hfil\cr
   \global\let\br\br@label #1\cr %upper labels
   \global\let\br\br@down #1\cr   %upper braces
   #2\cr % main line of the formula
   \global\let\br\br@up #3\cr %lower braces
   \global\let\br\br@label #3\cr   %lower labels
 }}% finished whole alignment and \vbox
 \advance\dimen@\ht\z@ %calc. the necessary lowering
 \lower\dimen@\hbox{\box\z@} % shift tho get the desired height
}}
\def\br@up#1#2{\multispan{#1}\upbracefill}
\def\br@down#1#2{\multispan{#1}\downbracefill}
\def\br@label#1#2{\multispan{#1}\hidewidth $#2$\hidewidth}
\makeatother

\begin{document}

%%%%%%%%%%% 
%% Title %%
%%%%%%%%%%%
\title[Decomposition of Generalized Theta Functions]
{A Canonical Decomposition of Generalized Theta Functions on
 the Moduli Stack of Gieseker Vector Bundles}
\author[Ivan Kausz]{Ivan Kausz}
\thanks{Partially supported by the DFG} 
\date{August 28, 2004}
\address{NWF I - Mathematik, Universit\"{a}t Regensburg, 93040 Regensburg, 
Germany}
\email{ivan.kausz@mathematik.uni-regensburg.de}
%\begin{abstract}
%\end{abstract}

\maketitle

\tableofcontents

%\setcounter{section}{-1}

%--------------------------------------------------------------------
\section{Introduction}

In this paper we prove a decomposition formula for generalized
theta functions which is motivated by what in conformal field
theory  is called the 
{\em factorization rule}.

A rational conformal field theory associates a finite dimensional 
vector space (called the {\em space of conformal blocks})
to a pointed nodal projective algebraic curve over $\C$
whose marked points are labeled by elements of a certain finite set.
If we choose a singular point $p$ in such a labeled pointed curve $X$,
then we get a new pointed curve $\tilde{X}$
by taking the partial normalization at $p$ and marking all points
which lie either over one of the marked points of $X$ (old points)
or over the singularity $p$ (two new points).
The factorization rule gives a canonical direct sum decomposition of the space 
of conformal blocks associated to $X$ with its labeled marked points,
such that the summands appearing in that decomposition are spaces
of conformal blocks associated to the pointed curve $\tilde{X}$ whose
old marked points are labeled by the same elements as the corresponding
points of $X$. The direct sum runs over a certain finite set of labeling of
the two new points.         

In the case of a Wess-Zumino-Witten conformal field theory associated
to a simply connected semi-simple algebraic group $G$ and a natural
number $\kappa\geq 1$ 
Tsuchiya, Ueno and Yamada have given a mathematical
definition of the spaces of conformal blocks in terms of the representation 
theory of the affine Lie algebra associated to $G$ and they have shown that
these spaces satisfy the factorization rule 
(\cite{TUY}, \cite{Ueno}, cf. also \cite{Sorger} for an overview).
It has been conjectured by physicists and later proved by 
various mathematicians that the spaces of conformal blocks
of Tsuchiya, Ueno, Yamada have an algebro-geometric interpretation:
They can be identified with spaces of global sections 
(called {\em generalized theta functions}) of certain
line bundles on the moduli space (or moduli stack) of $G$-bundles
with parabolic structures at the marked points 
(\cite{F1}, \cite{BL}, \cite{LS}, \cite{T}). 
It should be noted that in the case of {\em singular} curves the
moduli spaces of $G$-bundles is {\em non-compact}.
Nevertheless, for semi-simple $G$ the space of global sections
of a line bundle on these moduli spaces is still finite-dimensional. 
This has been pointed out to me by the referee and follows e.g. 
from \cite{T}, Theorem 3.

Generalized theta functions can also be defined on moduli of principal
$\GL_n$-bundles (or equivalently of vector bundles) on a smooth curve.
For singular curves however the theta line bundle over the non-compact
moduli space of vector bundles carries too many sections, so one has
to compactify it, to get a reasonable notion of generalized theta functions.
There exist at least two approaches to compactify the moduli space of vector
bundles (of given rank and degree, say) on a singular curve.
One construction uses torsion free sheaves 
(\cite{Seshadri1}, \cite{Newstead}, \cite{F2}),
the other one  works with 
certain vector bundles on modifications on the singular curve
and
has been introduced by Gieseker \cite{Gieseker}
in the rank two case and has been generalized to arbitrary rank
by Nagaraj and Seshadri \cite{NS} and myself \cite{degeneration}.
The torsion free sheaves approach works for arbitrary singularities;
Gieseker's approach has up to now been carried out
only for the case where the curve is irreducible with only one 
ordinary double point. 

A version of the factorization rule for generalized theta functions
has been formulated and proved 
by Narasimhan, Ramadas and Sun
in the framework of moduli varieties of semi-stable torsion free sheaves 
of fixed rank and degree
(\cite{NR}, \cite{Ramadas}, \cite{Sun1}, \cite{Sun2}).
They also prove 
(at least for rank$=2$ or genus$\geq 4$) that
in case of a one-dimensional family of curves which is generically
smooth and degenerates at one point to an irreducible nodal curve with
one singularity, the spaces of generalized theta functions form
the fibers of a finite rank vector bundle on the one-dimensional base.

In the present paper we prove a factorization rule 
for generalized theta functions on the moduli stack of Gieseker
vector bundles on an irreducible curve with one node 
(a stack version of Gieseker's approach which
I have constructed in \cite{degeneration}).
Our result is somewhat stronger than the analogous result of 
Narasimhan, Ramadas and Sun, since we obtain a {\em canonical}
decomposition, whereas the decomposition proved by those
authors is non-canonical.

In the last chapter we show that generalized theta functions
on the moduli stack of Gieseker vector bundles behave well
in families. However, since we are dealing with Artin-stacks
which are neither separated nor of finite type over the base,
we can not argue by cohomological flatness, but have to make
an explicit dimension calculation using the Verlinde formula
for $\SL_n$. So our result cannot be regarded as an alternative
to the representation theoretic approach to the Verlinde
factorization. Rather it shows that the Gieseker type 
``compactification'' of the moduli stack of vector bundles on
a singular curve leads to the ``correct selection rules'' for the
sections of the determinant line bundle.

Here is the main result 
of the paper (cf. Theorem \ref{main}):

\vspace{3mm}
\noindent
{\bf Theorem:}
{\em
Let $k$ be an algebraically closed field of characteristic zero.
Let $C_0$ be an irreducible projective 
algebraic curve over $k$ with one ordinary
double point, let $\tilde{C_0}$ be its normalization and let
$p_1, p_2$ be the two points of $\tilde{C_0}$ which are mapped
onto the singular point of $C_0$.
Let $\Theta$ be the theta line bundle on
$\GVB$, the moduli stack of rank $n$ Gieseker vector bundles on 
$C_0$ and let $\kappa$ be a positive integer.
Then there is a canonical isomorphism of $k$-vector spaces
$$
H^0(\GVB,\Tk)\Isomto
\Oplus_{(a,b)\in A'}
H^0(\PB,\Tk_{\PB}(a,b))
\quad.
$$
Here $\PB$ is the moduli stack parametrizing vector bundles  on 
$\tilde{C_0}$ together with full flags in the fibers at
the points $p_1$ and $p_2$ and $A'$ is a finite set (depending on $\kappa$) 
which parametrizes a set of line bundles $\Tk_{\PB}(a,b)$ on $\PB$.
}
\vspace{3mm}

The main ingredient of the proof of Theorem
\ref{main} is a result from my earlier
paper \cite{degeneration}, which says that $\GVB$ has normal
crossing singularities and that there is a diagram of algebraic
$k$-stacks:
$$
\xymatrix{
\VB &
\GVBD \ar[l]_f \ar[r]^\nu &
\GVB
}
$$
where $\VB$ is the moduli stack of rank $n$ 
vector bundles on $\tilde{C_0}$,
the morphism $f$ is a locally trivial fibration whose standard fiber
is a certain canonical compactification $\KGLn$ 
of the general linear group $\GL_n$ and
the morphism $\nu$ identifies the stack $\GVBD$ with
the normalization of $\GVB$.

In view of that diagram the strategy of the proof of the theorem
is quite straightforward:
We identify the space $H^0(\GVB,\Tk)$ with the subspace
of $H^0(\GVBD,\nu^*\Tk)$ consisting of sections of $\nu^*\Tk$ whose values
coincide at points which map onto the same point of $\GVB$.
We show that the line bundle $\nu^*\Theta$ is naturally
isomorphic to $f^*\Thetat\tensor\Delta$ where 
$\Thetat$ is the theta line bundle on $\VB$ and
$\Delta$ is a line bundle whose restriction to a fiber of $f$ is a fixed
line bundle on the compactification $\KGL_n$.
We then apply the result from \cite{cohomology} where we have
decomposed the cohomology of line bundles on $\KGL_n$ in terms
of irreducible representations of $\GL_n\times\GL_n$.
This yields a canonical decomposition 
$$
H^0(\GVBD,\nu^*\Tk)\Isomto
\Oplus_{(a,b)\in A(\Delta^k)}
H^0(\PB,\Tk_{\PB}(a,b))
\quad.
$$
Finally we determine how the subspace $H^0(\GVB,\Tk)$
of $H^0(\GVBD,\nu^*\Tk)$ behaves with respect to this decomposition.
It turns out that the composite morphism
$$
H^0(\GVB,\Tk)\injto
H^0(\GVBD,\nu^*\Tk)\Isomto
\!\!\!\!\!
\Oplus_{(a,b)\in A(\Delta^k)}
\!\!\!\!\!
H^0(\PB,\Tk_{\PB}(a,b))\To
\!\!\!
\Oplus_{(a,b)\in A'}
\!\!\!
H^0(\PB,\Tk_{\PB}(a,b))
$$
is an isomorphism. The last arrow in this diagram is simply
the projection induced by the inclusion of the finite sets
$A'\subset A(\Dk)$.

The greater part of this work has been carried out during a stay
at the Tata Institute of Fundamental Research in Bombay.
Its hospitality is gratefully acknowledged.
I am deeply indebted to Don Zagier, who provided me with an ingenious proof
of Lemma \ref{zagier}.

%weiter
%--------------------------------------------------------------------------
\section{Complements on modifications of pointed nodal curves}
\label{complements}

Let $S$ be an arbitrary scheme or algebraic stack.
In this section we present some constructions which yield two-pointed
nodal curves over $S$.

\begin{definition}
\begin{enumerate}
\item
A {\em nodal curve} over $S$ is a morphism
$\pi:\Cc\to S$ which is projective, finitely
presented and flat and whose geometric fibers are reduced curves
with only ordinary double points as singularities.
We require furthermore that
for each point $z\in S$ we have
$H^0(\Cc_z,\Oo_{\Cc_z})=\kappa(z)$,
where $\Cc_z$ denotes the fiber of $\pi$ at $z$ and $\kappa(z)$ 
is the residue field of the point $z$.
\item
A {\em one-pointed nodal curve} over $S$ is a tuple 
$(\Cc,\pi,s)$, where $\pi:\Cc\to S$ is a nodal curve
and where
$s$ is a section of $\pi$, whose image is contained in
the smooth locus of $\pi$
such that for each point $z\in S$ we have 
$
H^0(\Cc_z,\Oo_{\Cc_z}(-s(z)))=(0)
$.
The morphism $\pi$ will often be omitted from the notation.
\item
A {\em two-pointed nodal curve} over $S$ is a tuple 
$(s_1,\Cc,\pi,s_2)$, where $\pi:\Cc\to S$ is a nodal curve
and where
$s_1$ and $s_2$ are disjoint sections of $\pi$ 
such that
$(\Cc,\pi,s_1)$ and $(\Cc,\pi,s_2)$ are one-pointed nodal curves.
We will often write $(s_1,\Cc,s_2)$ instead
of $(s_1,\Cc,\pi,s_2)$.
\item
A morphism $(s_1,\Cc,\pi,s_2)\to(s'_1,\Cc',\pi',s'_2)$
of two-pointed nodal curves over $S$ is an $S$-morphism
$f:\Cc\to\Cc'$ with 
$s'_i=f\comp s_i$ for $i=1,2$.
\end{enumerate}
\end{definition}

\begin{definition}
Let $(s_1,\Cc,s_2)$ and $(t_1,\Dc,t_2)$ be two-pointed nodal curves over $S$.
Then we define the two-pointed nodal curve 
$$
(s_1,\Cc,s_2)\bot(t_1,\Dc,t_2):=
(r_1,\Bc,r_2)
\quad,
$$
where $\Bc$ is the curve $(\Cc\sqcup\Dc)/(s_2=t_1)$ and 
the sections $r_1$, $r_2$ are defined by
$
r_1: 
S \overset{s_1}{\to}
\Cc\setminus \{s_2(S)\}\injto
\Bc 
$
and
$
r_2: 
S \overset{t_2}{\to}
\Dc\setminus \{t_1(S)\} \injto
\Bc 
$ 
respectively.
\end{definition}

\begin{definition}
\label{|L,M|}
Let $L_1$ and $L_2$ be two line bundles on $S$.
We denote by $|L_1,L_2|$ the two-pointed nodal curve 
$(s_1,\Cc,\pi,s_2)$, where $\Cc:=\Pp(L_1\oplus L_2)$ and
$s_i:S\to \Cc$ is defined by the invertible quotient
$L_1\oplus L_2\to L_i$ for $i=1,2$.
\end{definition}

\begin{lemma}
\label{easy}
Let $L$ and $M$ be two line bundles on $S$ and let
$(s,\Cc,\pi,t):=|L,M|$. Let 
$\pi^*(L\oplus M)\to\Oo_{\Cc}(1)$ 
be the tautological
invertible quotient on $\Cc$. 
Then there are canonical isomorphisms
\begin{enumerate}
\item
$s^*\Oo_{\Cc}(s)=L\tensor M^{-1}$
\quad and\quad
$t^*\Oo_{\Cc}(t)=L^{-1}\tensor M$,
\item
$
\Oo_{\Cc}(1)=
\Oo_{\Cc}(s)\tensor\pi^*M=
\Oo_{\Cc}(t)\tensor\pi^*L
$.
\end{enumerate}
\end{lemma}

\begin{proof}
This follows easily from the universal property of $\Pp(L\oplus M)$.
\end{proof}

For the convenience of the reader I will now recall a definition and a
result from \cite{degeneration} \S 5 which are needed in the following.

\begin{definition}
\label{def mpc}
Let $S$ be a scheme and 
$(\Cc,\pi,s)$ a one-pointed nodal curve.
A {\em simple modification $(\Cc',f,\pi',s')$ of $(\Cc,\pi,s)$} 
is a diagram
$$
\xymatrix{
\Cc' \ar[rr]^f \ar[dr]^{\pi'} & & \Cc \ar[dl]_{\pi} \\
   & S \ar@/^/[ul]^{s'} \ar@/_/[ur]_s
}
$$
with the following properties:
\begin{enumerate}
\item
The triple $(\Cc',\pi',s')$ is again a one-pointed nodal curve over $S$.
\item
The diagram is commutative in the sense that $\pi\comp f=\pi'$
and $f\comp s'=s$.
\item 
The morphism $f$ is proper and finitely presented.
\item
Let $z\in S$ be a point. Then there are two possibilities for
the induced morphism $f_z:\Cc'_z\to\Cc_z$ of fibers over $z$:
Either $f_z$ is an isomorphism, or 
$\Cc'_z$ 
arises from
$R\isomorph\Pp_{\kappa(z)}^1$ and $\Cc_z$ by the identification
of a point in $R(\kappa(z))$ with $s(z)\in\Cc_z$,
and $f_z$ contracts $R$ to the point $s(z)$.
\end{enumerate}
\end{definition}

\begin{remark}
We will most often use a shorter expression by
saying that ``$(\Cc',s')$ is a simple modification
of $(\Cc,s)$'', the data $\pi'$, $\pi$ and $f$ being understood.
The definition implies, that $f$ induces an isomorphism 
$\Cc'\setminus f^{-1}(s(S))\isomto \Cc\setminus s(S)$
In particular, $\Cc\setminus s(S)$ can be considered as an open
subscheme of $\Cc'$.
\end{remark}

\begin{proposition}
\label{mod to Mmu}
Let $S$ be a scheme, and let $(\Cc,\pi,s)$ be a one-pointed
nodal curve over $S$.
Then there is a canonical 
isomorphism of groupoids:
$$
\left\{
\begin{array}{ll}
\text{Simple modifications of $(\Cc,\pi,s)$} \\
\text{in the sense of definition \ref{def mpc}}
\end{array}
\right\}
\isomto
\left\{
\begin{array}{llll}
\text{Pairs $( M,\mu)$, where $ M$ is } \\
\text{an invertible $\Oo_S$-module and}\\
\text{$\mu:\Oo_S\to  M$ is a global}\\
\text{section of $ M$}
\end{array}
\right\}
$$
\end{proposition}

\begin{sketch of proof}
\label{sketch}
For a proof I refer the reader to \cite{degeneration} \S 5.
Here  we only need the following details of the 
correspondence 
$
(C',\pi',f',s')\leftrightarrow(M,\mu)
$:

\begin{enumerate}
\item
Let $(\Cc',f,\pi',s')$ be a simple modification of $(\Cc,\pi,s)$.
Then $(M,\mu)$ is constructed as follows.
First of all we have $M=(s')^*\Oo_{\Cc'}(-s')\tensor s^*\Oo_{\Cc}(s)$.
For the section $\mu$
consider the following exact diagram
\[
\xymatrix@R=3ex{
0 \ar[r] &
\Oo_{\Cc'}(-s') \ar[r] \ar[d] &
\Oo_{\Cc'} \ar[r] \ar@{.>}[dl] \ar[d] &
s'_*s^*\Oo_{\Cc} \ar[r] \ar[d]^0 &
0 \\
0 \ar[r]&
\M \ar[r] &
f^*\Oo_{\Cc}(s) \ar[r] &
s'_*s^*\Oo_{\Cc}(s) \ar[r] & 
0 
}
\]
where
$\M:=\Oo_{\Cc'}(-s')\tensor f^*\Oo_{\Cc}(s)$ and where the vertical arrows
are induced by the morphism $\Oo_{\Cc'}\to f^*\Oo_{\Cc}(s)$ obtained by
applying the functor $f^*$ to the natural injection 
$\Oo_{\Cc}\injto\Oo_{\Cc}(s)$. 
Since the right vertical arrow obviously vanishes, the middle vertical arrow
factorizes as indicated by the dotted arrow. Applying $(s')^*$ th the morphism
$\Oo_{\Cc'}\to \M$ thus obtained, yields the section $\mu:\Oo_S\to M$.
\item
Let $(M,\mu)$ be an invertible $\Oo_S$-module, together with a section.
The associated nodal curve $\Cc'$ is canonically isomorphic to
$\Pp(\J)$, where $\J$ is the $\Oo_{\Cc}$-module
defined by the following exact sequence
\[
\xymatrix{
0\ar[r] &
\pi^*M^{-1} \ar[r]^(.3){(i,-\mu)} &
(\Oo_{\Cc}(s)\tensor\pi^*M^{-1})\oplus\Oo_{\Cc} \ar[r] &
\J \ar[r] &
0}
\]
where $i$ is the morphism induced by the inclusion
$\Oo_{\Cc}\injto\Oo_{\Cc}(s)$
and $-\mu$ is the negative of the morphism induced
by $\mu:\Oo_S\to M$.
\end{enumerate}
\end{sketch of proof}

\begin{definition}
\label{dashv}
Let $(s_1,\Cc,s_2)$ be a two-pointed nodal curve over $S$
and let $M$ be a line bundle on $S$ and $\mu$ a global section of $M$.
Then we denote by 
$$
(M,\mu)\dashv(s_1,\Cc,s_2)
$$
the two-pointed nodal
curve $(r_1,\Bc,r_2)$, where $(\Bc,r_1)$ is the simple modification
of $(\Cc,s_1)$ associated to the data $(M,\mu)$ 
by \ref{mod to Mmu}.
and $r_2$ is defined as the composition
$
r_2: S\overset{s_2}{\to}\Cc\setminus s_1(S) \injto \Bc 
$.
Similarly, we write 
$$
(s_1,\Cc,s_2)\vdash(M,\mu)
$$ 
for the two-pointed 
nodal curve $(t_1,\Dc,t_2)$, where $(\Dc,t_2)$ is the simple modification
of $(\Cc,s_2)$ associated to the data $(M,\mu)$ and $t_1$
is defined by the composition
$
t_1: S\overset{s_1}{\to}\Cc\setminus s_2(S) \injto \Dc 
$.
In situations where no doubts as to $\mu$ are likely to arise,
we will sometimes write 
$M\dashv(s_1,\Cc,s_2)$ and $(s_1,\Cc,s_2)\vdash M$
instead of 
$(M,\mu)\dashv(s_1,\Cc,s_2)$ and $(s_1,\Cc,s_2)\vdash(M,\mu)$.
\end{definition}

\begin{lemma}
\label{lemma bot}
Let $M$ be a line bundle on $S$ with zero section $0$  
and let $(s_1,\Cc,s_2)$ be
a two-pointed nodal curve. Then we have canonical isomorphisms
of two-pointed nodal curves as follows:
\begin{eqnarray*}
|\Oo_S,M\tensor s_1^*\Oo_{\Cc}(-s_1)|\bot(s_1,\Cc,s_2) &=& 
                                         (M,0)\dashv(s_1,\Cc,s_2) \\
(s_1,\Cc,s_2)\bot|M\tensor s_2^*\Oo_{\Cc}(-s_2),\Oo_S| &=& 
                                         (s_1,\Cc,s_2)\vdash(M,0)
\end{eqnarray*}
\end{lemma}

\begin{proof}
Let $M':=M\tensor s_1^*\Oo_{\Cc}(-s_1)$ and let
$(r_1,\Bc,r_2):=|\Oo_S,M'|\bot(s_1,\Cc,s_2)$.
Let $f:\Bc\to\Cc$ be the 
morphism, whose restriction to 
$\Pp(M'\oplus\Oo_S)$ 
is the structure morphism to $S$ composed with the section 
$s_1$ and whose restriction to $\Cc$ is
the identity morphism. Clearly $(f:\Bc\to\Cc,r_1)$ 
is a simple modification of $(\Cc,s_1)$. 
Furthermore, by \ref{easy} 
the $\Oo_S$-module $r_1^*\Oo_{\Bc}(-r_1)\tensor s_1^*\Oo_{\Cc}(s_1)$
is canonically isomorphic to $M$ and the canonical morphism
$\Oo_S\to r_1^*\Oo_{\Bc}(-r_1)\tensor s_1^*\Oo_{\Cc}(s_1)$ from construction
\ref{sketch} (1) vanishes. The canonical isomorphism 
$
|\Oo_S,M'|\bot(s_1,\Cc,s_2) 
= (M,0)\dashv(s_1,\Cc,s_2)
$
follows now from \ref{mod to Mmu}.
The other isomorphism follows completely analogously.
\end{proof}

\begin{lemma}
\label{switch}
Let $L$, $M$ be two line bundles on $S$ and let $\mu$ be a global
section of $M$. Then there is a canonical isomorphism of
two-pointed nodal curves as follows:
$$
|L,\Oo_S|\vdash (M,\mu)=
(M,\mu)\dashv|\Oo_S,(L\tensor M)^{-1}|
\quad.
$$
\end{lemma}

\begin{proof}
Let 
\begin{eqnarray*}
(s_1,\Cc_1,\pi_1,t_1) &:=& |L,\Oo_S| \\
(s,\Cc,\pi,t) &:=& |L,\Oo_S|\vdash (M,\mu) \\
(s_2,\Cc_2,\pi_2,t_2) &:=& |\Oo_S,(L\tensor M)^{-1}|
\end{eqnarray*}
and let $f_1: \Cc\to\Cc_1=\Pp(L\oplus\Oo_S)$ 
be the canonical morphism.
By \ref{sketch} (2)  we have $\Cc=\Pp(\J)$, where
$\J$ is the coherent $\Oo_{\Cc_1}$-module defined by the exact sequence
$$
\xymatrix{
0\ar[r] &
\pi_1^*M^{-1} \ar[r]^(.3){(i,-\mu)} &
(\Oo_{\Cc_1}(t_1)\tensor\pi_1^*M^{-1})\oplus\Oo_{\Cc_1} \ar[r]^(.8){p} &
\J \ar[r] &
0}
$$
Now consider the composite morphism
$$
\xymatrix{
\pi^*(L\tensor M)^{-1}\oplus\Oo_{\Cc}  \ar[r]^(.45){f_1^*j} &
(f_1^*\Oo_{\Cc_1}(t_1)\tensor\pi^* M^{-1})\oplus \Oo_{\Cc} 
                                        \ar[r]^(.75){f_1^*p} &
f_1^*\J \ar@{->>}[r]^(.47)q &
\Oo_{\Cc/\Cc_1}(1)
}
$$
where 
$
j:\pi_1^*(L\tensor M)^{-1}\oplus\Oo_{\Cc_1}
  \to (\Oo_{\Cc_1}(t_1)\tensor\pi_1^*M^{-1})\oplus\Oo_{\Cc_1}
$ 
is induced by the injection 
$
\pi_1^* L^{-1}=\pi_1^*t_1^*\Oo_{\Cc_1}(t_1)=\Oo_{\Cc_1}(t_1-s_1)
\injto\Oo_{\Cc_1}(t_1)
$
and $q$ is the tautological invertible quotient on $\Cc=\Pp(\J)$.
It is easy to check that the morphism $p\comp j$ is surjective.
Hence the morphism $u:=q\comp f_1^*p\comp f_1^* j$ is surjective 
and defines a morphism 
$$
\xymatrix{
f_2:\Cc\ar[r] &
\Pp(\Oo_S\oplus(L\tensor M)^{-1})=\Cc_2
}
\quad.
$$

It is not hard to see that the pull-back by $s$ and $t$ of the epimorphism
$u$ identifies with the epimorphism
$
\Oo_S\oplus(L\tensor M)^{-1}\to\Oo_S
$
(projection to the first component)
and 
$
\Oo_S\oplus(L\tensor M)^{-1}\to(L\tensor M)^{-1}
$
(projection to the second component)
respectively.
Therefore we have $s_2=f_2\comp s$ and $t_2=f_2\comp t$.

By considering the case where $S$ is the spectrum of a field one
verifies that $(\Cc,s)$ is a simple modification of $(\Cc_2,s_2)$.

Clearly we have $s^*\Oo_{\Cc}(-s)\tensor s_2^*\Oo_{\Cc_2}(s_2)=M$
and by going through the constructions it follows that 
the canonical morphism 
$\Oo_S\to s^*\Oo_{\Cc}(-s)\tensor s_2^*\Oo_{\Cc_2}(s_2)$
from \ref{sketch} (1) identifies with $\mu$.

The lemma now follows from \ref{mod to Mmu}.
\end{proof}

\begin{proposition}
\label{switches}
Let $M_0,\dots,M_q$ be invertible $\Oo_S$-modules and for $i\in[1,q]$
let $\mu_i$ be a global section of $M_i$. Then the two-pointed nodal
curve
$$
|M_0,\Oo_S|\vdash(M_1,\mu_1)\vdash(M_2,\mu_2)\vdash\dots\vdash
(M_q,\mu_q)
$$
is canonically isomorphic to the two-pointed nodal curve
$$
(M_1,\mu_1)\dashv(M_2,\mu_2)\dashv\dots\dashv
(M_q,\mu_q)\dashv|\Oo_S,\Tensor_{i=0}^q M_i^{-1}|
\quad.
$$
\end{proposition}

\begin{proof}
This follows by a $q$-fold application of lemma \ref{switch}.
\end{proof}

The next proposition is needed in the proof of proposition
\ref{isomorphism}.  

\begin{proposition}
\label{Sigma}
Let $(s_1,\Cc,s_2)$ be a two-pointed nodal curve such that $\Cc\to S$
is smooth and let $(L_i,\lambda_i)_{i=1,\dots,q}$ and 
$(M_i,\mu_i)_{i=1,\dots,r}$ be two families of invertible $\Oo_S$-modules
with sections.
Consider the two-pointed nodal curve
\[
(t_1,\Bc,t_2):=
(M_r,\mu_r)\dashv\dots\dashv(M_1,\mu_1)\dashv
(s_1,\Cc,s_2)
\vdash(L_1,\lambda_1)\vdash\dots\vdash(L_q,\lambda_q)
\]
and let $\Sigma\injto\Bc$ be the singular locus 
of the morphism $\Bc\to S$.
Then we can express $\Sigma$ as a disjoint union of closed subschemes
of $S$ as follows: 
\[
\Sigma=
\left(
\bigsqcup_{i=1}^{q}\{\lambda_i=0\}
\right)
\sqcup
\left(
\bigsqcup_{i=1}^{r}\{\mu_i=0\}
\right)
\quad.
\]
\end{proposition}

\begin{proof}
The proposition follows from repeated application of the following
assertion:

Let $(\Cc',\pi',f,s')$ be a simple modification of the one-pointed
nodal curve $(\Cc,\pi,s)$ over $S$ and let $(M,\mu)$ be the corresponding
line bundle with section. Let $\Sigma\subset\Cc$ and
$\Sigma'\subset\Cc'$ be the
locus of non-smoothness of $\pi$ and $\pi'$ respectively.
Then we have
\[
\Sigma'=\Sigma \sqcup Y
\quad,
\]
where $Y\subseteq S$ is the closed subscheme of $S$
defined by the equation $\mu=0$.

Since $\Sigma$ and the section $s(S)$ are disjoint closed subschemes of 
$\Cc$, and $\Cc'\setminus f^{-1}(s(S))$ is isomorphic to $\Cc\setminus s(S))$,
we can identify $\Sigma$ with the locus of non-smoothness of 
$\Cc'\setminus f^{-1}(s(S))\to S$ and
it suffices to show that the locus of non-smoothness of 
$\Cc'\setminus\Sigma\to S$ is isomorphic to $Y$.
For this we may replace $\Cc$ by an open affine neighborhood
$V=\Spec(B)\subseteq\Cc\setminus\Sigma$
of the section $s(S)$ and $S$ by $U=s^{-1}(V)$. Then $U$ is
also affine, say $U=\Spec(A)$ and we may also assume that
$M|_U$ is trivial and that $\mu$ is given by $a\in A$.
By \cite{degeneration}, 5.5 we have
$f^{-1}(V)=\Proj(R)$ for the graded $B$-algebra
$$
R:=B[X,Y]/(bX-aY)
\quad,
$$
where $b\in B$ is the regular element generating
the ideal associated to the closed subscheme $s(U)\subset V$.
Now $W:=\Proj(R)$ is the union of the open affine pieces
$W_1=\Spec(R_1)$ and $W_2=\Spec(R_2)$, where
$$
R_1\isomorph B[T]/(b-aT)
\qquad
\text{and}
\qquad
R_2\isomorph B[T]/(bT-a)
\quad.
$$
A simple calculation shows that the $R_1$-module 
$\Omega^1_{R_1/A}$  is trivial 
of rank one and that the first Fitting ideal of the $R_2$-module
$\Omega^1_{R_2/A}$ is $(b,T)$.
Therefore $\Sigma'\cap W\isomorph\Spec(R_2/(b,T))=\Spec(A/a)=Y\cap U$.
\end{proof}

%----------------------------------------------------------------------------
\section{Review of a compactification of $\GL_n$}
\label{review kgl}

Let $k$ be a field and $E$, $F$ two $k$-vector spaces of rank $n$.
In \cite{kgl} I have studied a certain compactification
$\KGL(E,F)$ of the scheme $\Isom(E,F)$ of isomorphisms from $E$ to $F$
which has properties similar to the so called wonderful compactification
of adjoint linear groups introduced by De Concini and Procesi.
In particular, the complement of $\Isom(E,F)$ in $\KGL(E,F)$ is a divisor
with normal crossings whose irreducible components are smooth.
As in \cite{kgl}, I denote these components by 
$Y_0,\dots,Y_{n-1},Z_0,\dots,Z_{n-1}$.

More generally, the construction of $\KGL(E,F)$ works also
in a relative situation (cf. \cite{kgl} \S 9). 
Thus let $S$ be a scheme and let
$E$ and $F$ be two locally free $\Oo_S$-modules of rank $n$.
Then there is a natural $S$-scheme $\KGL(E,F)$ containing the scheme
$\Isom(E,F)$ as an open subscheme such that for each point $z\in S$
the fiber of $\KGL(E,F)\to S$ over $z$ is the compactification
$\KGL(E_z,F_z)$ of $\Isom(E_z,F_z)$. We denote by 
$Y_0,\dots,Y_{n-1},Z_0,\dots,Z_{n-1}$ the divisors in $\KGL(E,F)$ which
are fiber-wise the components of the boundary of $\Isom(E_z,F_z)$ in
$\KGL(E_z,F_z)$.
The construction may be extended to the case where $S$ is an algebraic
stack.

The main theorem in \cite{kgl} is a concrete description of the $T$-valued
points of $\KGL(E,F)$ (for any $S$-scheme $T$). As we need this result
in the sequel, I will recall the necessary definitions.

Let $T$ be an $S$-scheme and let $\E$, $\F$ be two locally free
$\Oo_T$-modules of rank $n$. A 
{\em bf-morphism} from $\E$ to $\F$ is a tuple
$
g=(L, \lambda, \E\to\F, \F\to M\tensor\E, r)
$
where $L$ is an invertible $\Oo_T$-module, $\lambda$ is a section of $L$,
the arrows $\E\to\F$ and $\F\to L\tensor\E$ are $\Oo_T$-module morphisms and 
$r$ is an integer between $0$ and $n$ 
such that locally on $T$ there exist isomorphisms
\begin{eqnarray*}
\E &\isomto& r\Oo_T\oplus(n-r)\Oo_T \\
\F &\isomto& r\Oo_T\oplus(n-r)L
\end{eqnarray*}
with the property that via these isomorphisms the morphisms
$
\E\to\F
$
and 
$\F\to L\tensor\E$
are expressed by the diagonal matrices
$$
\text{
$
\left[
\begin{array}{cc}
\eins_r & 0 \\
0 & \lambda\eins_{n-r}
\end{array}
\right]
$
\quad and \quad
$
\left[
\begin{array}{cc}
\lambda\eins_r & 0 \\
0 & \eins_{n-r}
\end{array}
\right]
$
}
$$
respectively.
We will often use the following more suggestive notation for the 
bf-morphism $g$:
$$
g=\left(
\xymatrix@C2.5ex{
\E
\ar[rr]^r_{(L,\lambda)}
& & 
\F
\ar@/_1.2pc/|{\tensor}[ll]
}
\right)
\quad.
$$

Let $T$, $\E$, $\F$ be as above.
A {\em generalized isomorphism} from $\E$ to $\F$ is a sequence of
bf-morphisms connected as follows:
$$
\xymatrix@C=1.2ex{
\E
\ar@/^1.2pc/|{\tensor}[rr]
& &
E_1 
\ar[ll]_0^{(M_0,\mu_0)}
\ar@/^1.2pc/|{\tensor}[rr]
& &
E_2
\ar[ll]_1^{(M_1,\mu_1)}
& 
\dots
& 
E_{n-1}
\ar@/^1.2pc/|{\tensor}[rrr]
& & &
E_n
\ar[lll]_{n-1}^{(M_{n-1},\mu_{n-1})}
\ar[r]^\sim
& 
F_n
\ar[rrr]^{n-1}_{(L_{n-1},\lambda_{n-1})}
& & &
F_{n-1}
\ar@/_1.2pc/|{\tensor}[lll]
& 
\dots
&
F_2
\ar[rr]^1_{(L_1,\lambda_1)}
& &
F_1
\ar[rr]^0_{(L_0,\lambda_0)}
\ar@/_1.2pc/|{\tensor}[ll]
& &
\F
\ar@/_1.2pc/|{\tensor}[ll]
}
$$
which has properties for which we refer the reader to
\cite{kgl} 5.2, since they will not be of importance here.

Now let $S$, $E$, $F$ be as in the beginning of this section 
and let $T$ be an $S$-scheme.
The main theorem in \cite{kgl} is the following:

\begin{theorem}
There is a natural bijection
between the set of 
$T$-valued points of $\KGL(E,F)$ and the set of (equivalence classes of)
generalized isomorphisms from $E_T$ to $F_T$, where $E_T$ and $F_T$
denote the pull back of $E$ and $F$ to $T$.
\end{theorem}

In particular, if $f:\KGL(E,F)\to S$ denotes the structure morphism,
then there exists a universal generalized isomorphism
$$
\xymatrix@C=1.2ex{
f^*E
\ar@/^1.2pc/|{\tensor}[rr]
& &
E_1 
\ar[ll]_0^{(M_0,\mu_0)}
\ar@/^1.2pc/|{\tensor}[rr]
& &
E_2
\ar[ll]_1^{(M_1,\mu_1)}
& 
\dots
& 
E_{n-1}
\ar@/^1.2pc/|{\tensor}[rrr]
& & &
E_n
\ar[lll]_{n-1}^{(M_{n-1},\mu_{n-1})}
\ar[r]^\sim
& 
F_n
\ar[rrr]^{n-1}_{(L_{n-1},\lambda_{n-1})}
& & &
F_{n-1}
\ar@/_1.2pc/|{\tensor}[lll]
& 
\dots
&
F_2
\ar[rr]^1_{(L_1,\lambda_1)}
& &
F_1
\ar[rr]^0_{(L_0,\lambda_0)}
\ar@/_1.2pc/|{\tensor}[ll]
& &
f^*F
\ar@/_1.2pc/|{\tensor}[ll]
}
$$
from $f^* E$ to $f^* F$.

For each pair of subsets $I,J\subseteq[0,n-1]$ let 
$\Ob_{I,J}=\Ob_{I,J}(E,F)$ be the closed subscheme of
$\KGL(E,F)$ defined by the equations
$\mu_i=0$ $(i\in I)$ and $\lambda_j=0$ $(j\in J)$.
It is non-empty if and only if $\min(I)+\min(J)\geq n$.
With this notation we have
$$
\Ob_{I,J}=
\left(\bigcap_{i\in I} Z_i\right)
\cap
\left(\bigcap_{j\in J} Y_j\right)
\quad.
$$
In particular, we have 
$\KGL(E,F)=\Ob_{\emptyset,\emptyset}$, 
$Y_i=\Ob_{\emptyset,\{i\}}$ 
and
$Z_i=\Ob_{\{i\},\emptyset}$.

%----------------------------------------------------------------------------
\section{Review of the cohomology of line bundles on $\KGL(E,F)$}
\label{cohomology}

Let $k$ be a field of characteristic zero.
Throughout this section, $S$ will denote a $k$-scheme (or more generally
an algebraic $k$-stack).
We fix two locally free $\Oo_S$-modules $E$ and $F$ of rank $n$.
Let $\KGL(E,F)$ be the compactification of $\Isom(E,F)$ introduced
in \S \ref{review kgl} and denote by $f:\KGL(E,F)\to S$ the structure
morphism. Let
$$
\xymatrix@C=1.2ex{
f^*E
\ar@/^1.2pc/|{\tensor}[rr]
& &
E_1 
\ar[ll]_0^{(M_0,\mu_0)}
\ar@/^1.2pc/|{\tensor}[rr]
& &
E_2
\ar[ll]_1^{(M_1,\mu_1)}
& 
\dots
& 
E_{n-1}
\ar@/^1.2pc/|{\tensor}[rrr]
& & &
E_n
\ar[lll]_{n-1}^{(M_{n-1},\mu_{n-1})}
\ar[r]^\sim
& 
F_n
\ar[rrr]^{n-1}_{(L_{n-1},\lambda_{n-1})}
& & &
F_{n-1}
\ar@/_1.2pc/|{\tensor}[lll]
& 
\dots
&
F_2
\ar[rr]^1_{(L_1,\lambda_1)}
& &
F_1
\ar[rr]^0_{(L_0,\lambda_0)}
\ar@/_1.2pc/|{\tensor}[ll]
& &
f^*F
\ar@/_1.2pc/|{\tensor}[ll]
}
$$
be the universal generalized isomorphism from $f^*E$ to
$f^*F$.

\begin{lemma}
\label{detE_n}
We have the following canonical isomorphisms
of invertible $\Oo_{\KGL(E,F)}$-modules:
\begin{eqnarray*}
f^*(\det E)\tensor\Tensor_{i=0}^{n-1}M_i^{i-n}\  =\ 
\det E_n
&=&
\det F_n \ =\  
f^*(\det F)\tensor\Tensor_{i=0}^{n-1}L_i^{i-n}
\end{eqnarray*}
\end{lemma}

\begin{proof}
For $i\in[0,n-1]$ we denote by $g_i$ and $h_i$ the
bf-morphism
$$
\xymatrix@C=2.5ex{
E_{i+1}
\ar[rr]^i_{(M_i,\mu_i)}
& &
E_i
\ar@/_1.2pc/|{\tensor}[ll]
}
\qquad\text{and}\qquad
\xymatrix@C=2.5ex{
F_{i+1}
\ar[rr]^i_{(L_i,\lambda_i)}
& &
F_i
\ar@/_1.2pc/|{\tensor}[ll]
}
$$
respectively.
By Proposition 6.2 in \cite{kgl} these induce canonical morphisms
\begin{eqnarray*}
\wedge^{-n}g_i &:& \det(E_i) \To M_i^{n-i}\tensor\det(E_{i+1}) \\
\wedge^{n}h_i &:& \det(F_{i+1}) \To L_i^{i-n}\tensor\det(E_{i}) 
\end{eqnarray*}
respectively. 
It follows that we have canonical morphisms
\begin{eqnarray*}
g :=
(\wedge^{-n}g_{n-1})\comp
(\wedge^{-n}g_{n-2})\comp
\dots\comp
(\wedge^{-n}g_{0}) &:&
f^*\det(E) 
\To
\det(E_n)\tensor\Tensor_{i=0}^{n-1}M_i^{n-i}
\\
h :=
(\wedge^{n}h_{0})\comp
(\wedge^{n}h_{1})\comp
\dots\comp
(\wedge^{n}h_{n-1}) &:&
\det(F_n)
\To
f^*\det(F)\tensor\Tensor_{i=0}^{n-1}L_i^{i-n}
\end{eqnarray*}
Let $\varphi:\det(E_n)\isomto\det(F_n)$ be the isomorphism
induced by the isomorphism $E_n\isomto F_n$.

By \cite{kgl} 6.5 the morphism 
$$
\wedge^n\Phi=h\comp\varphi\comp g: 
f^*\det(E)\To\Tensor_{i=0}^{n-1}(M_i^{n-i}\tensor L_i^{i-n})
$$
is nowhere vanishing and consequently $g$ and $h$ are isomorphisms.
\end{proof}

Recall from \S \ref{review kgl} that for each pair of subsets
$I,J\subseteq[0,n-1]$
with $\min(I)+\min(J)\geq n$ the closed subscheme 
$\Ob_{I,J}$ of $\KGL(E,F)$ is defined as the 
zero locus of the sections $\mu_i$ ($i\in I$) and
$\lambda_j$ ($j\in J$).
Let
$$
\xymatrix@R=1.1ex{
i_{I,J}: \Ob_{I,J}\ar@{^{(}->}[r] &
\KGL(E,F) \\
f_{I,J}: \Ob_{I,J}\ar[r] &
S
}
$$
denote the inclusion morphisms and structure morphisms respectively.

Let 
$\Fl:=\text{Fl}(E)\times_S\text{Fl}(F)$,
where $\text{Fl}(E)$ and $\text{Fl}(F)$ denote the varieties over $S$ which
parametrize full flags in $E$ and $F$ respectively.
Let 
$
0=\E_0\subset\E_1\subset\dots\subset\E_n=E\tensor\Oo_{\Fl}
$
and
$
0=\F_0\subset\F_1\subset\dots\subset\F_n=F\tensor\Oo_{\Fl} 
$
be the two universal flags on $\Fl$.
For 
$a=(a_1,\dots,a_n)\in \Z^n$,\ 
$b=(b_1,\dots,b_n)\in \Z^n$
we define the invertible $\Oo_{\Fl}$-module
$$
\Oo_{\Fl}(a,b) := \Tensor_{i=1}^n(\E_i/\E_{i-1})^{\tensor a_i} \tensor
                    \Tensor_{i=1}^n(\F_i/\F_{i-1})^{\tensor b_i} 
\quad,
$$
which we sometimes abbreviate by $\Oo(a,b)$ if no confusion
is likely to arise.

\begin{definition}
\label{A}
Let $L$ be a line bundle on $\KGL(E,F)$ of the form
$$
L=\Tensor_{i=0}^{n-1}(M_i^{m_i}\tensor L_i^{l_i})\tensor
  f^*(\det E)^e\tensor f^*(\det F)^d
\quad.
$$ 
Let $I,J\subseteq[0,n-1]$ and let $i_1:=\min(I)$, $j_1:=\min(J)$
where it is understood that $\min(\emptyset)=n$. Assume $i_1+j_1\geq n$.
We denote by $A_{IJ}(L)$ the set of all elements 
$(a,b)\in \Z^n\times\Z^n$, which have
the following properties:
\begin{enumerate}
\item
$a_1\leq a_2\leq\dots\leq a_n$ 
\item
$
\sum_{j=i+1}^{n}(a_j-e)\leq m_i
$
for all $i\in [n-j_1,n-1]$
and equality holds for $i\in I$.
\item
$
\sum_{j=1}^{n-i}(a_j-e)\geq -l_i
$
for all $i\in[n-i_1,n-1]$
and equality holds for $i\in J$.
\item
For all $i\in[1,n]$ the 
equality $a_i-e=-b_{n-i+1}+d$ holds.
\end{enumerate}
For abbreviation we will often write $A(L)$ instead of 
$A_{\emptyset,\emptyset}(L)$.
\end{definition}

\begin{theorem}
\label{bundle decomposition}
Let $L$ be a line bundle on $\KGL(E,F)$ of the form
$$
L=\Tensor_{i=0}^{n-1}(M_i^{m_i}\tensor L_i^{l_i})\tensor
  f^*(\det E)^e\tensor f^*(\det F)^d
$$
and let $I,J\subseteq[0,n-1]$ be subsets with $\min(I)+\min(J)\geq n$.
Then the following holds:

1.
The $\Oo_S$-module $(f_{IJ})_*i_{IJ}^*L$ is locally free and
comes with a canonical decomposition as follows:
$$
(f_{IJ})_*i_{IJ}^*L=
\Oplus_{(a,b)\in A_{IJ}(L)}(f_{\Fl})_*\Oo_{\Fl}(a,b)
\quad,
$$
where $f_{\Fl}:\Fl\to S$ denotes the structure morphism.

2.
\label{compatibility}
The decomposition stated in 1.
is compatible with restriction in the sense
that the following diagram commutes:
$$
\xymatrix{
f_*L \ar@{=}[d] \ar[rr]^{\text{Res}} & &
(f_{I,J})_*i_{I,J}^*L \ar@{=}[d] \\
\text{$
\underset{(a,b)\in A(L)}{\Oplus}(f_{\Fl})_*\Oo_{\Fl}(a,b)
$} 
\ar@<1ex>@{->>}[r] &
\text{$
\underset{(a,b)\in A(L)\cap A_{I,J}(L)}{\Oplus} (f_{\Fl})_*\Oo_{\Fl}(a,b)
$} 
\ar@<1ex>@{^(->}[r] &
\text{$
\underset{(a,b)\in A_{I,J}(L)}{\Oplus} (f_{\Fl})_*\Oo_{\Fl}(a,b)
$}
}
$$
where the lower arrows are the canonical projection and inclusion morphisms
induced by the inclusions $A(L)\cap A_{I,J}(L)\subseteq A(L)$ and
$A(L)\cap A_{I,J}(L)\subseteq A_{I,J}(L)$ respectively. 

3. Let 
$$
L'=\Tensor_{i=0}^{n-1}(M_i^{m'_i}\tensor L_i^{l'_i})\tensor
  f^*(\det E)^e\tensor f^*(\det F)^d
\quad,
$$
where $m'_i\leq m_i$ and $l'_j\leq l_j$ and equality holds,
if $i\in I$ and $j\in J$ respectively.
The following diagram commutes:
$$
\xymatrix{
(f_{I,J})_*i_{I,J}^*L' 
\ar@{^(->}[r]^{\tensor\mu^{m-m'}\tensor\lambda^{l-l'}} 
\ar@{=}[d]
&
(f_{I,J})_*i_{I,J}^*L 
\ar@{=}[d]
\\
\text{$
\underset{(a,b)\in A_{I,J}(L')}{\Oplus}(f_{\Fl})_*\Oo_{\Fl}(a,b)
$}
\ar@<1ex>@{^(->}[r]
&
\text{$
\underset{(a,b)\in A_{I,J}(L)}{\Oplus}(f_{\Fl})_*\Oo_{\Fl}(a,b)
$}
}
$$
where the upper horizontal arrow is induced by the section
$$
\left.
\left(
\mu_0^{m_0-m'_0}\tensor\dots\tensor
\mu_{n-1}^{m_{n-1}-m'_{n-1}}\tensor
\lambda_0^{l_0-l'_0}\tensor\dots\tensor
\lambda_{n-1}^{l_{n-1}-l'_{n-1}}
\right)\right|_{\Ob_{I,J}}
$$
of $i_{I,J}^*(L\tensor (L')^{-1})$
and the lower horizontal arrow is induced by the inclusion
$A_{I,J}(L')\subseteq A_{I,J}(L)$.
\end{theorem}

\begin{proof}
This is an easy consequence of the main result in \cite{cohomology}.
\end{proof}

%---------------------------------------------------------------------------
\section{Review of moduli of Gieseker vector bundles}
\label{review degeneration}

Let $k$ be an algebraically closed field of characteristic zero.
Let $C_0$ be an irreducible projective curve over $k$
which is smooth except for one ordinary double point $p\in C_0(k)$.
Let $\Ct$ be the normalization of $C_0$ and 
let $p_1,p_2$ be the two $k$-valued points of $\Ct$ lying above the
singular point $p$.

For an integer $q\geq 1$ let $C_q$ be the curve which arises from
$C_0$ by inserting a chain of length $q$ of projective lines at
the point $p$:

\vspace{3mm}
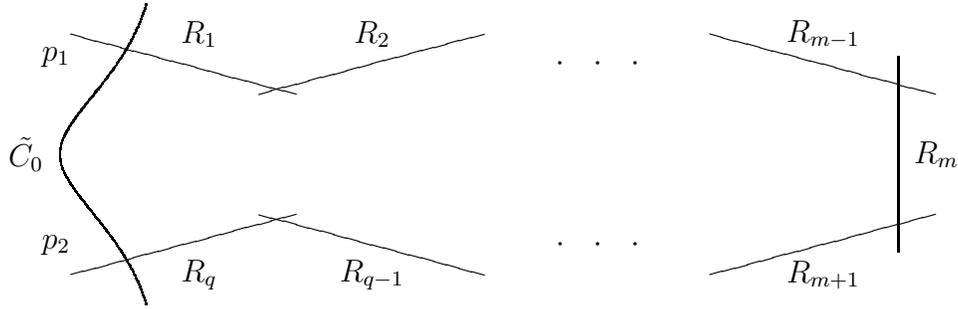
\begin{center}
\parbox{8cm}{
\xy <0mm,-20mm>; <1mm,-20mm>:
(10,20);(10,-20)**\crv{(7,10)&(-10,0)&(7,-10)};
(0,16);(30,8) **@{-}; (17,16) *{\text{$R_1$}};
(25,8);(55,16) **@{-}; (40,16) *{\text{$R_2$}};
(65,12)*{\cdot};
(70,12)*{\cdot};
(75,12)*{\cdot};
(85,16);(115,8) **@{-}; (100,16) *{\text{$R_{m-1}$}};
(110,13);(110,-13) **@{-}; (115,0) *{\text{$R_m$}};
(0,-16);(30,-8) **@{-}; (17,-16) *{\text{$R_q$}};
(25,-8);(55,-16) **@{-}; (40,-16) *{\text{$R_{q-1}$}};
(65,-12)*{\cdot};
(70,-12)*{\cdot};
(75,-12)*{\cdot};
(85,-16);(115,-8) **@{-}; (100,-16) *{\text{$R_{m+1}$}};
(-2,13) *{\text{$p_1$}};
(-2,-12) *{\text{$p_2$}};
(-6,0) *{\text{$\Ct$}};
\endxy
}
\end{center}
\vspace{3mm}

\begin{definition}
1.
A {\em Gieseker vector bundle of rank $n$ on $C_0$} is a
pair $(X\to C_0,\F)$, where $X=C_r$ for some $r\in[0,n]$ and
$\F$ is a vector bundle of rank $n$ on $X$ such that the following holds:
\begin{enumerate}
\item
The morphism $X\to C_0$ is the identity if $r=0$, and it 
contracts the chain of projective lines
into the singular point $p$ of $C_0$ if $r\geq 1$.
\item
The restriction of $\F$ to any of the inserted projective lines
$R_i$ is of the form
$$
d_i\Oo_{R_i}(1)\oplus(n-d_i)\Oo_{R_i}
$$
for some $d_i\geq 1$.
\item
If $r\geq 1$
let $R=\bigcup_iR_i$ be the inserted chain of projective lines
and denote by $p_1,p_2\in R$ the points at which it meets the
curve $\Ct$. Then we have
$$
H^0(R,\F|_R(-p_1-p_2))=(0)
\quad.
$$
\end{enumerate}

2.
Let $T$ be a $k$-scheme. A {\em Gieseker vector bundle of rank $n$
on $C_0$ over $T$} is a pair $(X\to C_0\times T,\F)$, 
where $X\to C_0\times T$ is a morphism of curves
over $T$ and $\F$ is a vector bundle on $X$ 
such that if $z$ is a point in $T$ and if we denote by
$X_z$ the fiber of $X\to T$ at 
$z$, then the pair 
$(X_z\to C_0\tensor_k\kappa(z), \F|_{X_z})$ is a 
Gieseker vector bundle on $C_0\tensor_k\kappa(z)$.
We will often write $(X,\F)$ instead of $(X\to\C_0\times T,\F)$.
\end{definition}

\begin{definition}
A {\em Gieseker vector bundle data of rank $n$ on $C_0$ 
over a $k$-scheme $T$} is a triple 
$(X\to C_0\times T,\F,x)$, where
$(X\to C_0\times T,\F)$ is a Gieseker vector bundle data
or rank $n$ on $C_0$ and $x:T\to X$ is a section of $X\to T$
whose image is in the singular locus of $X\to T$.
We will often write $(X,\F,x)$ instead of $(X\to C_0\times T,\F,x)$.
\end{definition}

\begin{remark}
Let
$(X\to C_0\times T,\F,x)$ 
be a Gieseker vector bundle data over $T$.
Then there is a canonical two-pointed curve 
$(x_1,\Cc,x_2)$ over $T$ such that $X$ can be
constructed from $\Cc$ by identifying the two sections
$x_1$ and $x_2$. Indeed, the curve $\Cc$ is simply the
blow-up of $X$ along the closed subscheme $x(T)$ and 
$x_1(T)\sqcup x_2(T)$ is the pre-image of $x(T)$.
The composition $\Cc\to X\to C_0\times T$ factorizes naturally
through a morphism $\Cc\to \Ct\times T$ and
the pull back $\F'$ of $\F$ to $\Cc$ comes with a natural
isomorphism $\varphi:x_1^*\F'\isomto x_2^*\F'$.
The datum 
$((x_1,\Cc,x_2)\to(p_1,\Ct\times T,p_2),\F',\varphi)$
is equivalent to the datum
$(X\to C_0\times T,\F,x)$.

\end{remark}

In \cite{degeneration}
I have shown that there are algebraic moduli stacks $\GVB$ and $\GVBD$ 
parametrizing Gieseker vector bundles and Gieseker vector bundle data
of rank $n$ on $C_0$ respectively. Furthermore, the
stack $\GVBD$ is smooth, the stack $\GVB$ has normal crossing 
singularities and the forgetful morphism
$$
\nu':
\left\{
\begin{array}{ll}
\GVBD &\to \GVB \\
(X,\F,x) &
\mapsto(X,\F)
\end{array}
\right.
$$ 
identifies the stack $\GVBD$
with the normalization of $\GVB$.

Let $\VB$ be the moduli stack of vector bundles of rank $n$ on $\Ct$
and denote by $\pit:\Ct\times\VB\to\VB$ the projection onto the
second factor.
Let $E$ and $F$ be the pull back of the universal 
vector bundle on $\Ct\times\VB$ via the section of $\pit$
induced by the point $p_1$ and $p_2$ respectively.
The main result in \cite{degeneration} is the construction of a 
canonical isomorphism 
$$
\tau: \KGL:=\KGL(E,F)\isomto\GVBD
\quad.
$$
The isomorphism $\tau$ is defined by a certain family 
$$
(\Cf'\to\KGL,\ \E',\ s:\KGL\to\Cf' )
$$ 
of Gieseker vector bundle data on $\KGL$.
In the remainder of this section I will recall some details of 
the construction of this family. 

Let $f:\KGL\to\VB$ be the structure morphism and let
$\Phi:$
$$
\xymatrix@C=1.2ex{
f^*E
\ar@/^1.2pc/|{\tensor}[rr]
& &
E_1 
\ar[ll]_0^{(M_0,\mu_0)}
\ar@/^1.2pc/|{\tensor}[rr]
& &
E_2
\ar[ll]_1^{(M_1,\mu_1)}
& 
\dots
& 
E_{n-1}
\ar@/^1.2pc/|{\tensor}[rrr]
& & &
E_n
\ar[lll]_{n-1}^{(M_{n-1},\mu_{n-1})}
\ar[r]^\sim
& 
F_n
\ar[rrr]^{n-1}_{(L_{n-1},\lambda_{n-1})}
& & &
F_{n-1}
\ar@/_1.2pc/|{\tensor}[lll]
& 
\dots
&
F_2
\ar[rr]^1_{(L_1,\lambda_1)}
& &
F_1
\ar[rr]^0_{(L_0,\lambda_0)}
\ar@/_1.2pc/|{\tensor}[ll]
& &
f^*F
\ar@/_1.2pc/|{\tensor}[ll]
}
$$
be the universal generalized isomorphism from $f^*E$ to $f^*F$.
Let $(s_1,\Bf,s_2)$ be the two-pointed nodal curve
(with the notation of \S \ref{complements})
\[
(M_{n-1},\mu_{n-1})\dashv\dots\dashv(M_0,\mu_0)\dashv
(p_1,\Ct\times\KGL,p_2)
\vdash(L_0,\mu_0)\vdash\dots\vdash(L_{n-1},\lambda_{n-1})
\quad.
\]
Then $\Cf'$ is the curve $\Bf/(s_1=s_2)$ over $\KGL$ and 
$s:\KGL\to \Cf$ is the composition 
$$
\xymatrix{
\KGL \ar[r]^{s_i} &
\Bf \ar[r] &
\Cf' 
}
\quad
\text{(for $i=1$ or $i=2$).}
$$

Let $\Et$ be the universal vector bundle on $\Ct\times\VB$
and let ${\Et}':=(\id_{\Ct}\times f)^*\Et$ be its pull-back
to $\Ct\times\KGL$.
The generalized isomorphism $\Phi$ together with the 
vector bundle ${\Et}'$ induce a vector bundle
$\G$ of rank $n$ on the two-pointed nodal curve 
$(s_1,\Bf,s_2)$. The details of the construction of $\G$
out of $\Phi$ and ${\Et}'$ are given in \S 7 and \S 9 of 
\cite{degeneration}. For the purpose of this paper 
it suffices to know that $\G$ has the property that
$$
{\Et}'=(h_*\G(-s_1-s_2))(p_1+p_2)
\quad,
$$
where $h:\Bf\to\Ct\times\KGL$ is the canonical projection.
Furthermore there are canonical isomorphisms
$$
s_1^*\G=E_n
\qquad
\text{and}
\qquad
s_2^*\G=F_n
\quad.
$$

The vector bundle $\E'$ on $\Cf'$ is constructed from $\G$ by
using the isomorphism $E_n\isomto F_n$, which is part of the
data contained in $\Phi$, to glue together $\G$ along the sections
$s_1$ and $s_2$.

For future reference we collect the various curves, bundles, 
and moduli spaces
in the following diagram:
$$
\xymatrix@R=1.5ex{
\text{$\Et$} &&
\text{${\Et}'$}  & &
\G &&
\E'&&
\E
\\
\text{$\Ct$}\times \VB \ar[dd]^{\pit} & &
\text{$\Ct$}\times \KGL \ar[ll] \ar[dd]^{\pit'} & &
\Bf \ar[ll]_h \ar[dd]^{\rho} \ar[rr]^g & &
\Cf' \ar[dd]_{\pi'} \ar[rr] &&
\Cf  \ar[dd]_{\pi}
\\ \\
\VB  \ar@/^/[uu]^{p_1} \ar@<-1ex>@/_/[uu]_{p_2} &&
\KGL \ar[ll]_f \ar@{=}[rr] \ar@/^/[uu]^{p_1} \ar@<-1ex>@/_/[uu]_{p_2} &&
\KGL \ar@{=}[rr] \ar@/^/[uu]^{s_1} \ar@<-1ex>@/_/[uu]_{s_2} &&
\KGL \ar[rr]^{\nu} \ar@/_/[uu]_{s} &&
\GVB
\\
E,F & &
f^*E, f^*F & &
E_n, F_n & &
s^*\E'
}
$$
Here we have set $\nu:=\nu'\comp\tau$. The pair 
$(\Cf,\E)$ is the universal Gieseker vector bundle on $C_0$ 
over $\GVB$. The outer squares in this diagram are Cartesian.
The bundles $\E'$ and $\G$ are the pull back of the bundle
$\E$ and the bundle ${\Et}'$ is the pull back of the bundle
$\Et$. The bundles on the bottom are the pull back by the respective
sections of the bundles upstairs.

%---------------------------------------------------------------------------
\section{A remarkable set of isomorphisms}
\label{remarkable}

We keep the notation from \S \ref{review degeneration}.
Thus we have a diagram
$$
\xymatrix{
\VB &
\KGL \ar[l]_f \ar[r]^{\nu} &
\GVB
}
$$
of algebraic stacks, where $\GVB$ is the moduli stack of Gieseker vector
bundles of rank $n$ on the singular curve $C_0$, and 
$\VB$ is the moduli stack of vector bundles of rank $n$ on the
normalization $\Ct$ of $C_0$.
We have $\KGL=\KGL(E,F)$, where $E$ and $F$ are the pull-back
of the universal bundle on $\Ct\times\VB$ along the sections $p_1$ and 
$p_2$ respectively.
Let 
$$
\xymatrix@C=1.2ex{
f^*E
\ar@/^1.2pc/|{\tensor}[rr]
& &
E_1 
\ar[ll]_0^{(M_0,\mu_0)}
\ar@/^1.2pc/|{\tensor}[rr]
& &
E_2
\ar[ll]_1^{(M_1,\mu_1)}
& 
\dots
& 
E_{n-1}
\ar@/^1.2pc/|{\tensor}[rrr]
& & &
E_n
\ar[lll]_{n-1}^{(M_{n-1},\mu_{n-1})}
\ar[r]^\sim
& 
F_n
\ar[rrr]^{n-1}_{(L_{n-1},\lambda_{n-1})}
& & &
F_{n-1}
\ar@/_1.2pc/|{\tensor}[lll]
& 
\dots
&
F_2
\ar[rr]^1_{(L_1,\lambda_1)}
& &
F_1
\ar[rr]^0_{(L_0,\lambda_0)}
\ar@/_1.2pc/|{\tensor}[ll]
& &
f^*F
\ar@/_1.2pc/|{\tensor}[ll]
}
$$
be the universal generalized isomorphism 
from $f^* E$ to $f^* F$.
As in \S \ref{review kgl} let 
$Y_i$ and $Z_i$ be the divisors in $\KGL$ defined by 
the equations $\lambda_i=0$ and $\mu_i=0$ respectively
and let 
\begin{eqnarray*}
i_{Y_i}:=i_{\emptyset,\{i\}} &:& Y_i\injto\KGL \\
i_{Z_i}:=i_{\{i\},\emptyset} &:& Z_i\injto\KGL \\
\end{eqnarray*}
be the respective inclusion morphisms.
In the following proposition we compute the fibre product
$\KGL\times_{\GVB}\KGL$. It is the key ingredient in proposition
\ref{H^0 as subspace}, where we identify the space of global
sections of a line bundle $\Ll$ on $\GVB$ with a subspace 
of the space of global sections of its pull-back to $\KGL$.  
\begin{proposition}
\label{isomorphism}
After the choice of isomorphisms $\m_{\Ct,p_i}/\m_{\Ct,p_i}^2=k$
for $i=1,2$ we have:

1.
For each $j\in[0,n-1]$ there is a canonical isomorphism
$\beta_j: Y_j\isomto Z_j$,
which makes the following diagram
commutative:
$$
\xymatrix{
Y_j
\ar[rr]^{\beta_j} \ar[dr]_{\nu_{Y_j}} 
& &
Z_j \ar[dl]^{\nu_{Z_j}} 
\\
& 
\GVB 
&
}
\quad,
$$
where $\nu_{Y_j}:=\nu\comp i_{Y_j}$ and $\nu_{Z_j}:=\nu\comp i_{Z_j}$.

2. \label{pullback via beta}
Let $i,j\in[0,n-1]$. Then we have
\begin{eqnarray*}
\beta_j^*i_{Z_j}^*(M_i,\mu_i)
&=&
\left\{
\begin{array}{lll}
i_{Y_j}^* 
(M_{n-j+i},\mu_{n-j+i}) & \text{for $i\in[0,j-1]$} 
\\
i_{Y_j}^* 
(\Tensor_{r=0}^{n-1}(M_r\tensor L_r)^{-1},0) & \text{for $i=j$} 
\\
i_{Y_j}^* 
(L_{n+j-i},\lambda_{n+j-i}) & \text{for $i\in[j+1,n-1]$}
\end{array}
\right.
\\
\beta_j^*i_{Z_j}^*(L_i,\lambda_i)
&=&
\left\{
\begin{array}{ll}
(\Oo,1) & \text{for $i\in[0, n-j-1]$} 
\\
i_{Y_j}^* 
(L_{j+i-n},\lambda_{j+i-n}) & \text{for $i\in[n-j,n-1]$}
\end{array}
\right.
\end{eqnarray*}

3.
The following morphism is an isomorphism:
$$
(\id,\id)\sqcup
\bigsqcup_{i=0}^{n-1}(i_{Y_i},i_{Z_i}\comp\beta_i)
\sqcup
\bigsqcup_{i=0}^{n-1}(i_{Z_i},i_{Y_i}\comp\beta^{-1}_i)
\ :\ 
\KGL\sqcup
\bigsqcup_{i=0}^{n-1}Y_i
\sqcup
\bigsqcup_{i=0}^{n-1}Z_i
\To
\KGL\times_{GVB}\KGL
$$
\end{proposition}

\begin{proof}
Let $T$ be a $k$-scheme. 
A $T$-valued point of $Y_j$ is given by a pair $(\F,\Psi)$, where
$\F$ is a vector bundle of rank $n$
on $\Ct\times S$ and $\Psi$ is a generalized isomorphism
$$
\xymatrix@C=1.2ex{
G_0
\ar@/^1.2pc/|{\tensor}[rr]
& &
G_1 
\ar[ll]_0^{(A_0,a_0)}
\ar@/^1.2pc/|{\tensor}[rr]
& &
G_2
\ar[ll]_1^{(A_1,a_1)}
& 
\dots
& 
G_{n-1}
\ar@/^1.2pc/|{\tensor}[rrr]
& & &
G_n
\ar[lll]_{n-1}^{(A_{n-1},a_{n-1})}
\ar[r]^\sim
& 
H_n
\ar[rrr]^{n-1}_{(B_{n-1},b_{n-1})}
& & &
H_{n-1}
\ar@/_1.2pc/|{\tensor}[lll]
& 
\dots
&
H_2
\ar[rr]^1_{(B_1,b_1)}
& &
H_1
\ar[rr]^0_{(B_0,b_0)}
\ar@/_1.2pc/|{\tensor}[ll]
& &
H_0
\ar@/_1.2pc/|{\tensor}[ll]
}
$$
from $G_0:=p_1^*\F$ to $H_0:=p_2^*\F$,
such that $b_j=0$.

Let $(X\to C_0\times T,\Hc,x)$ be the Gieseker vector
bundle data associated to $(\F,\Psi)$ by the canonical
isomorphism $\tau:\KGL\to\GVBD$.
Recall from  \S \ref{review degeneration} that $X$ is constructed
from the two-pointed nodal curve
$(x_1,\Cc,x_2)=$ 
$$
(A_{n-1},a_{n-1})\dashv\dots\dashv
(A_0,a_0)\dashv
(p_1,\Ct\times T,p_2)\vdash\
(B_0,b_0)\vdash\dots\vdash\
(B_{n-1},b_{n-1})
$$
by identifying the two sections $x_1$ and $x_2$.

We define the two-pointed nodal curves
\begin{eqnarray*}
(r_1,\Bc,r_2) &:=&
A_{n-1}\dashv\dots\dashv
A_{n-j}\dashv
(p_1,\Ct\times S,p_2)\vdash B_0\vdash\dots\vdash
B_{j-1} \\
(t_1,\Dc,t_2) &:=&
|\Tensor_{r=0}^{j}B_r,\Oo_S|\vdash
B_{j+1}\vdash\dots\vdash
B_{n-1}
\end{eqnarray*}

Since the sections $a_0,\dots,a_{n-j-1}$ are nowhere vanishing, and
since we have $b_j=0$
and $r_2^*\Oo_{\Bc}(-r_2)=\Tensor_{r=0}^{j-1}B_r$
(cf. \ref{sketch}, here we make use of the identification 
$\m_{\Ct,p_2}/\m_{\Ct,p_2}^2=k$),
it follows from \ref{lemma bot} that we have
$$
(x_1,\Cc,x_2)=(r_1,\Bc,r_2)\bot(t_1,\Dc,t_2)
\quad.
$$

Now we define
$$
(x'_1,\Cc',x'_2):=(t_1,\Dc,t_2)\bot(r_1,\Bc,r_2)
\quad.
$$
Then the curve $X'$ obtained from $\Cc'$ by identifying the sections
$x'_1$ and $x'_2$ is canonically isomorphic to $X$.
Thus we have a new Gieseker vector bundle data 
$(X'\to C_0\times T,\Hc',x')$,
where $X'=X$, $\Hc'=\Hc$ and $x'$ is the composition
$$
\xymatrix{
T \ar[r]^{x'_m} &
\Cc' \ar[r] &
X'=X
}
\qquad
\text{($m=1$ or $m=2$).}
$$

Via the isomorphism $\tau:\KGL\isomto\GVBD$ there corresponds 
to $(X',\Hc',x')$ a
$T$-valued point of $\KGL$ which is given by a pair $(\F',\Psi')$,
where $\F'$ is a vector bundle on $\Ct$ and $\Psi'$ is a 
generalized isomorphism from $p_1^*\F'$ to $p_2^*\F'$.

Recall from \S \ref{review degeneration} that
if we write $\Psi'=$
$$
\xymatrix@C=1.2ex{
G'_0
\ar@/^1.2pc/|{\tensor}[rr]
& &
G'_1 
\ar[ll]_0^{(A'_0,a'_0)}
\ar@/^1.2pc/|{\tensor}[rr]
& &
G'_2
\ar[ll]_1^{(A'_1,a'_1)}
& 
\dots
& 
G'_{n-1}
\ar@/^1.2pc/|{\tensor}[rrr]
& & &
G'_n
\ar[lll]_{n-1}^{(A'_{n-1},a'_{n-1})}
\ar[r]^\sim
& 
H'_n
\ar[rrr]^{n-1}_{(B'_{n-1},b'_{n-1})}
& & &
H'_{n-1}
\ar@/_1.2pc/|{\tensor}[lll]
& 
\dots
&
H'_2
\ar[rr]^1_{(B'_1,b'_1)}
& &
H'_1
\ar[rr]^0_{(B'_0,b'_0)}
\ar@/_1.2pc/|{\tensor}[ll]
& &
H'_0
\ar@/_1.2pc/|{\tensor}[ll]
}
$$
then 
$
(x'_1,\Cc',x'_2)
$
is isomorphic to
$$
(A'_{n-1},a'_{n-1})\dashv\dots\dashv
(A'_0,a'_0)\dashv
(p_1,\Ct\times S,p_2)\vdash
(B'_0,b'_0)\vdash\dots\vdash
(B'_{n-1},b'_{n-1})
\quad.
$$
Since by \ref{switches} we have 
$$
(t_1,\Dc,t_2) =
B_{j+1}\dashv\dots\dashv
B_{n-1}\dashv
|\Oo_S,\Tensor_{r=0}^{n-1}B_r^{-1}|
\quad,
$$
it follows from \ref{lemma bot} that 
\begin{eqnarray}
\label{A'}
(A'_i,a'_i):=
\left\{
\begin{array}{lll}
(A_{n-j+i},a_{n-j+i})  & 
\text{for $i\in[0,j-1]$} \\
(\Tensor_{r=0}^{n-1}(A_r\tensor B_r)^{-1},0) &
\text{for $i=j$} \\
(B_{n+j-i},b_{n+j-i}) &
\text{for $i\in[j+1,n-1]$}
\end{array}
\right.
\end{eqnarray}
(for $i=j$ we have made use of the identification
$\m_{\Ct,p_1}/\m_{\Ct,p_1}^2=k$ and the fact that 
the $a_0,\dots,a_{n-j-1}$ are nowhere vanishing)
and
\begin{eqnarray}
\label{B'}
(B'_i,b'_i):=
\left\{
\begin{array}{ll}
(\Oo_S,1)  & 
\text{for $i\in[0,n-j-1]$} \\
(B_{j-n+i},b_{j-n+i}) &
\text{for $i\in[n-j,n-1]$}
\quad.
\end{array}
\right.
\end{eqnarray}

In particular, we have $a'_j=0$ and therefore $(\F',\Psi')$ is
in fact a $T$-valued point of the closed substack $Z_j$ of $\KGL$.
We define $\beta_j: Y_j\to Z_j$ by the rule
$$
(\F,\Psi)\mapsto(\F',\Psi')
\quad.
$$
Since the inverse of $\beta_j$ can be constructed completely analogously,
it is clear that $\beta_j$ is an isomorphism.
By construction, we have 
$$
\nu_{Y_j}(\F,\Psi)
=
(X\to C_0\times T,\F)
=
(X'\to C_0\times T,\F')
=
\nu_{Z_j}(\F',\Psi')=\nu_{Z_j}\comp\beta_j(\F,\Psi)
\quad.
$$
This shows the first part of the Proposition.

The second part follows from equations (\ref{A'}) and (\ref{B'}) above.

For the third part 
it is clearly sufficient to show that there exists a commutative
diagram of stacks:
\begin{eqnarray}
\label{triangle}
\xymatrix@R=1.5ex{
& &
\KGL\sqcup
\left(\bigsqcup_{i=0}^{n-1}Y_i\right)\sqcup
\left(\bigsqcup_{i=0}^{n-1}Z_i\right)
\ar[dd]^{
\id\sqcup
\left(\bigsqcup_{i=0}^{n-1}\beta_i\right)\sqcup
\left(\bigsqcup_{i=0}^{n-1}\beta_i^{-1}\right)}
\\
\KGL\times_{\GVB}\KGL
\ar[rru]^{\alpha}
\ar[rrd]^{\alpha'}
& &
\\
& &
\KGL\sqcup
\left(\bigsqcup_{i=0}^{n-1}Z_i\right)\sqcup
\left(\bigsqcup_{i=0}^{n-1}Y_i\right)
}
\end{eqnarray}
where the arrows $\alpha$ and $\alpha'$ are isomorphisms
such that the following holds:
\begin{eqnarray}
\label{alpha}
\pr_1\comp\alpha^{-1} &=&
\id\sqcup
\left(\bigsqcup_{i=0}^{n-1}i_{Y_i}\right)\sqcup
\left(\bigsqcup_{i=0}^{n-1}i_{Z_i}\right)
\\
\label{alpha'}
\pr_2\comp(\alpha')^{-1} &=&
\id\sqcup
\left(\bigsqcup_{i=0}^{n-1}i_{Z_i}\right)\sqcup
\left(\bigsqcup_{i=0}^{n-1}i_{Y_i}\right)
\end{eqnarray}
Here, $\pr_m:\KGL\times_{\GVB}\KGL\to\KGL$ 
denotes the projection onto the $m$-the
factor $(m=1,2)$.

First we define the isomorphism $\alpha$.
Let $T$ be a scheme.
A $T$-valued point of $\KGL\times_{\GVB}\KGL$
is a pair $(\xi,\xi')$ of $T$-valued points of $\KGL$
such that $\nu(\xi)\isomorph(X,\Hc)\isomorph\nu(\xi')$
for some Gieseker vector bundle $(X,\Hc)$ over $T$.
The datum $(\xi,\xi')$ is equivalent
to the datum $(X,\Hc,x,x')$, where $x,x':T\to X$ 
are two sections of $X\to T$ whose image is contained
in the singular locus of $X\to T$ such that 
$(X,\Hc,x)=\tau(\xi)$ and $(X,\Hc,x')=\tau(\xi')$.
There are two cases:
\begin{description}
\item[First case]
If $x=x'$, then $(X,\Hc,x,x')$ is nothing else but a $T$-valued
point of $\GVBD\isomorph\KGL$.
\item[Second case]
If $x\neq x'$, then let $(x_1,\Cc,x_2)$ be the two-pointed nodal curve
over $T$
which is the partial normalization of $X$ along $x$.
The datum of $(X,\Hc,x,x')$ is clearly equivalent to the datum
$(\xi,T\to\Cc)$, where $T\to\Cc$ is a section of $\Cc\to T$
whose image is in the singular locus of $\Cc\to T$ such that
the composition $T\to\Cc\to X$ is the section $x'$.
But the datum $(\xi,T\to\Cc)$ describes precisely a $T$-valued point of
the closed substack $\Sigma$ of the curve $\Bf$ which is the locus of 
non-smoothness of the morphism $\Bf\to\KGL$.
From the definition of $\Bf$ and from \ref{Sigma} it follows
that we have 
$$
\Sigma=\left(\bigsqcup_{i=0}^{n-1}Z_i\right)\sqcup
\left(\bigsqcup_{i=0}^{n-1}Y_i\right)
$$
Thus in this case $(X,\Hc,x,x')$ is equivalent to a $T$-valued point of
the disjoint union of the $Y_i$ and $Z_i$.
\end{description}
Thus we have established an equivalence 
between $T$-valued points of the product $\KGL\times_{\GVB}\KGL$ 
and $T$-valued points  of the disjoint union of the stacks 
$\KGL$, $Y_i$ and $Z_i$.
This defines the isomorphism $\alpha$.
It is clear from the construction that equation (\ref{alpha})
holds.

The isomorphism $\alpha'$ is constructed similarly, with the
only difference that in the case $x\neq x'$ we take the
partial normalization $(x'_1,\Cc'x'_2)$ of $X$ along $x'$
and then define a $T$-valued point of $\Sigma$ by 
the datum $(\xi',T\to \Cc')$ where $T\to\Cc'$ is induced by $x$.
Again the equation (\ref{alpha'}) is clear.

The commutativity of the diagram (\ref{triangle})
is clear from the construction of the $\beta_i$.
\end{proof}

%----------------------------------------------------------------------------
\section{Decomposition of generalized theta functions}
\label{section decomposition}

We keep the notation from the end of \S \ref{review degeneration}.
In particular, $(\Cf,\E)$ denotes the universal 
Gieseker vector bundle over $\GVB$ and
$\pi:\Cf\to\GVB$ is the projection onto the base. 
$\Et$ is the universal vector bundle on $\Ct\times\VB$ and 
$\pit:\Ct\times\VB\to\VB$ is the projection onto the second factor.
Let $\Theta:=\det R\pi_*\E$ 
and $\tilde{\Theta}:=\det R\tilde{\pi}_*(\Et)$ be
the theta line bundle on $\GVB$ and on $\VB$ respectively.
Our convention for the determinant of the cohomology is
such that for a curve $X$ and a vector bundle $\F$ on $X$
we have $\det H(X,\F)=(\det H^0(X,\F))^{-1}\tensor\det H^1(X,\F)$.

We fix a positive integer $\kappa$.
Our aim is to decompose the space
$$
H^0(\GVB,\Tk)
$$
canonically into a direct sum, where the summands are related to
$\Thetat$.
The following proposition tells us that we can regard
$H^0(\GVB,\Tk)$ as a subspace of $H^0(\KGL,\nu^*(\Tk))$.

\begin{proposition}
\label{H^0 as subspace}
Let $\Ll$ be a line bundle on $\GVB$.
Then the canonical homomorphism 
$$
H^0(\GVB,\Ll)\to H^0(\KGL,\nu^*\Ll)
$$
is injective. A global section $\theta\in H^0(\KGL,\nu^*\Ll)$
is in the image of this homomorphism, if and only if
for each $j\in[0,n-1]$ the equality
$$
\beta_j^{\bullet}
i_{Z_j}^{\bullet}(\theta)=
i_{Y_j}^{\bullet}(\theta)
$$
holds,
where 
$i_{Y_j}^{\bullet}$ 
and 
$i_{Z_j}^{\bullet}$ 
denote the restriction homomorphisms
\begin{eqnarray*}
i_{Y_j}^{\bullet} &:& H^0(\KGL, \nu^*(\Ll))\to 
           H^0(Y_j, \nu_{Y_j}^*(\Ll))
\\
i_{Z_j}^{\bullet} &:& H^0(\KGL, \nu^*(\Ll))\to 
           H^0(Z_j, \nu_{Z_j}^*(\Ll))
\end{eqnarray*}
induced by $i_{Y_j}$ and $i_{Z_j}$ respectively,
and $\beta_j^{\bullet}$ denotes the pull-back isomorphism
$$
\beta_j^{\bullet}:H^0(Z_j, \nu_{Z_j}^*(\Ll)) \to
              H^0(Y_j, \nu_{Y_j}^*(\Ll))
$$
induced by $\beta_j$.
\end{proposition}

We need the following lemma, which is probably well-known, but for which
I did not find a reference.

\begin{lemma}
\label{descent}
Let $k$ be a field,
let $X$ be a smooth $k$-scheme and let $X_0\subset X$ be 
a divisor with normal crossings. Let $X_1\to X_0$ be the normalization
of $X_0$ and let $\Ll_0$ be an invertible $\Oo_{X_0}$-module.
Then the following sequence is exact:
$$
\xymatrix{
0 \ar[r] &
\text{$H^0(X_0,\Ll_0)$} 
\ar[r] &
\text{$H^0(X_1,\Ll_1)$} 
\ar@<.5ex>[r]  \ar@<-.5ex>[r] &
\text{$H^0(X_2,\Ll_2)$} 
}
$$
Here $X_2$ denotes the fiber product 
$X_1\times_{X_0}X_1$ and $\Ll_i$ is the pull back of
$\Ll_0$ by the morphism $X_i\to X_0$ for $i=1,2$.
The arrows are the obvious ones. 
\end{lemma}

\begin{proof}
{\em Step 1:} Assume $X=\Spec(R)$, where $R$ is a regular local ring and
$X_0=\Spec(R_0)$ where $R_0= R/(\prod_{i=1}^rx_i))$, the elements
$x_1,\dots,x_m$ form a regular system
of parameters for $R$ and $r\in[1,m]$.

By \cite{EGA} II 6.3.8 we have $X_1=\Spec(R_1)$ where 
$R_1=\prod_{i=1}^rR/(x_i)$ 
and it follows that $X_2=\Spec(R_2)$ where
$R_2=\prod_{i,j=1}^{r}R/(x_i,x_j)$.
We have to show the exactness of the sequence
$$
\xymatrix{
0 \ar[r] &
\text{$R_0$} \ar[r] &
\text{$R_1$} \ar@<.5ex>[r]  \ar@<-.5ex>[r] &
\text{$R_2$}
}
\quad.
$$

Since a regular local ring is a unique factorization domain,
any element of $R$ which is divisible by all $x_i$ for 
$i\in[1,r]$ is also divisible by the product $\prod_{i=1}^rx_i$.
This implies the injectivity of $R_0\to R_1$.

Let $(f_i)_{i\in[1,r]}$ be a family of elements in $R$ with
$f_i\equiv f_j$ mod $(x_i,x_j)$ for $i,j\in[1,r]$.
We have to show that there exists an element $f$ of $R$
with $f\equiv f_i$ mod $(x_i)$ for all $i\in[1,r]$.
If $r=1$, this statement is trivial, so assume that $r>1$ and
that there exists $f'\in R$ with 
$f'\equiv f_i$ mod $(x_i)$  for $i\in[1,r-1]$.
By assumption we have $f'-f_r\equiv g_ix_i$ mod $(x_r)$ 
for $i\in[1,r-1]$ and suitable $g_i\in R$. But the ring
$R/(x_r)$ is regular local and thus a unique factorization
domain and $x_1,\dots,x_{r-1}$ represent prime elements in $R/(x_r)$.
Therefore it follows that 
$f'-f_r= g\prod_{i=1}^{r-1}x_i-hx_r$
for suitable $g,h\in R$.
The element $f:=f'-g\prod_{i=1}^{r-1}x_i=f_r-hx_r$
has the required property.

\vspace{1mm}
{\em Step 2:} Assume that $X_0\subset X$ is a divisor with strict normal
crossings.

This means (cf. \cite{SGA I} XIII 2.1) that there exists a family
$(f_i)_{i\in I}$ of global sections of $\Oo_{X}$ indexed by a finite
set $I$, such that $X_0=\divi(\prod_{i\in I}f_i)$ and such that for
every $x\in X_0$ the closed subscheme of $X$ cut out by the ideal
$(f_i)_{i\in I(x)}$ is smooth of codimension equal to the cardinality 
of $I(x):=\{i\in I\ |\ f_i(x)=0\}$.

Let $\pi_i$ denote the morphism $X_i\to X_0 $ ($i=1,2$). By the first
step we have an exact sequence of $\Oo_{X_0}$-modules:
$$
\xymatrix{
0 \ar[r] &
\text{$\Oo_{X_0}$} \ar[r] &
\text{$(\pi_1)_*\Oo_{X_1}$} \ar@<.5ex>[r]  \ar@<-.5ex>[r] &
\text{$(\pi_2)_*\Oo_{X_2}$}
}
\quad.
$$
Tensoring with $\Ll_0$ 
and taking global sections yields the desired result.

\vspace{1mm}
{\em Step 3:} General case.

There exists an etale covering $X'\to X$ such that the pull back 
$X'_0\subset X'$ of $X_0$ is a divisor with strict normal crossings.
Let $X'':=X'\times_XX'$, denote by $X'_1$ and $X''_1$
the normalization of $X'_0$ and $X''_0$ respectively 
(which may be identified with the fiber product
 $X_1\times_{X_0}X'_0$ and $X_1\times_{X_0}X''_0$ respectively)
and let 
$X'_2:=X'_1\times_{X'_0}X'_1$, 
$X''_2:=X''_1\times_{X''_0}X''_1$.

Then we have the following commutative diagram:
$$
\vcenter{
\xymatrix@R=2ex{
& 0 \ar[d] & 0 \ar[d] & 0 \ar[d] 
\\
0 \ar[r] &
H^0(X_0) \ar[r] \ar[d] &
H^0(X_1) \ar@<.5ex>[r] \ar@<-.5ex>[r] \ar[d] &
H^0(X_2) \ar[d] 
\\
0 \ar[r] &
H^0(X'_0) \ar[r] \ar@<.5ex>[d] \ar@<-.5ex>[d] &
H^0(X'_1) \ar@<.5ex>[r] \ar@<-.5ex>[r] \ar@<.5ex>[d] \ar@<-.5ex>[d] &
H^0(X'_2) \ar@<.5ex>[d] \ar@<-.5ex>[d] 
\\
0 \ar[r] &
H^0(X''_0) \ar[r] &
H^0(X''_1) \ar@<.5ex>[r] \ar@<-.5ex>[r] &
H^0(X''_2) 
}
}
\eqno(*)
$$
where for a $X_0$-scheme $Y$ we denote by $H^0(Y)$ the space of
global sections of the pull back of $\Ll_0$ to $Y$.

Since the morphisms $X'_i\to X_i$ are etale coverings 
for $i=0,1,2$, the columns
of this diagram are exact. Since $X'_0$ and $X''_0$ are divisors 
of strict normal crossings in $X'$ and $X'\times_XX'$ respectively,
the second and third row in this diagram is also exact by step 2.
The exactness of the first row follows from this.
\end{proof}

\begin{proof}
(of the proposition \ref{H^0 as subspace}).

For abbreviation we set $X_0:=\GVB$ and $X_1:=\KGL$.
Let $X'_0\to X_0$ be a presentation of $\GVB$, i.e. $X'_0$ is a scheme and
$X'_0\to X_0$ is a smooth surjective morphism.
Let $X''_0:=X'_0\times_{X_0}X'_0$ and $X_2:=X_1\times_{X_0}X_1$.
Let 
$X'_i:=X_i\times_{X_0}X'_0$ and
$X''_i:=X_i\times_{X_0}X''_0$ for $i=1,2$.
Thus we have a diagram as follows:
$$
\xymatrix@R=2ex{
\text{$X''_2$} \ar@<.5ex>[r] \ar@<-.5ex>[r] \ar@<.5ex>[d] \ar@<-.5ex>[d] &
\text{$X''_1$} \ar[r]  \ar@<.5ex>[d] \ar@<-.5ex>[d] &
\text{$X''_0$} \ar@<.5ex>[d] \ar@<-.5ex>[d] 
\\
\text{$X'_2$} \ar@<.5ex>[r] \ar@<-.5ex>[r] \ar[d] &
\text{$X'_1$} \ar[r]  \ar[d] &
\text{$X'_0$} \ar[d] 
\\
\text{$X_2$} \ar@<.5ex>[r] \ar@<-.5ex>[r] &
\text{$X_1$} \ar[r] &
\text{$X_0$} 
}
$$
Taking $H^0$ of the pull back of $\Ll$ to the objects involved
in this diagram yields a diagram like $(*)$ in step 3 in the proof 
of lemma \ref{descent}.

Since $X'_i\to X_i$ is a presentation of the stack $X_i$ and
$X''_i=X'_i\times_{X_i}X'_i$, it follows from \cite{LM} 12.6.2 
that the columns in this diagram are exact.

Since $X'_0$ and $X''_0$ are normal crossing divisors in a smooth
 $k$-scheme (cf. \cite{degeneration} 3.21), since
$X'_1$ and $X''_1$ is the normalization of $X'_0$ and $X''_0$ 
respectively, and since we have 
$X'_2=X'_1\times_{X'_0}X'_1$ and
$X''_2=X''_1\times_{X''_0}X''_1$,
it follows from lemma \ref{descent} that the second and third row 
in $(*)$ are exact.

Therefore the first row in $(*)$, i.e. the sequence
$$
\xymatrix{
0 \ar[r] &
H^0(\GVB,\Ll) \ar[r] &
H^0(\KGL,\nu^*\Ll) \ar@<.5ex>[r] \ar@<-.5ex>[r] &
H^0(\KGL\times_{\GVB}\KGL,\Ll_2) 
}
\quad,
$$
where $\Ll_2$ is the pull back of $\Ll$ to 
$X_2=\KGL\times_{GVB}\KGL$, is also exact.

Proposition \ref{H^0 as subspace} now follows from 
Proposition \ref{isomorphism}.3.
\end{proof}

We will now study the space $H^0(\KGL,\nu^*(\Tk))$.
As a first step we compute the line bundle $\nu^*(\Theta)$.
The result is as follows:

\begin{proposition}
\label{pullback of Theta}
We have a canonical isomorphism of line bundles on $\KGL(E,F)$:
$$
\nu^*(\Theta)=\Delta\tensor f^*\Thetat
\quad,
$$
where $f$ is the morphism $\KGL\to\VB$ and where
$$
\Delta:=
(\Tensor_{i=0}^{n-1}M_i^{n-i})\tensor f^*(\det F)=
(\Tensor_{i=0}^{n-1}L_i^{n-i})\tensor f^*(\det E)
\quad.
$$
\end{proposition}

\begin{proof}
Recall from the end of \S \ref{review degeneration}
that we have a diagram of curves over $\VB$, $\KGL$ and $\GVB$ together
with vector bundles as follows:
$$
\xymatrix@R=1.5ex{
\text{$\Et$} &&
\text{${\Et}'$}  & &
\G &&
\E'&&
\E
\\
\text{$\Ct$}\times \VB \ar[dd]^{\pit} & &
\text{$\Ct$}\times \KGL \ar[ll] \ar[dd]^{\pit'} & &
\Bf \ar[ll]_h \ar[dd]^{\rho} \ar[rr]^g & &
\Cf' \ar[dd]_{\pi'} \ar[rr] &&
\Cf  \ar[dd]_{\pi}
\\ \\
\VB  \ar@/^/[uu]^{p_1} \ar@<-1ex>@/_/[uu]_{p_2} &&
\KGL \ar[ll]_f \ar@{=}[rr] \ar@/^/[uu]^{p_1} \ar@<-1ex>@/_/[uu]_{p_2} &&
\KGL \ar@{=}[rr] \ar@/^/[uu]^{s_1} \ar@<-1ex>@/_/[uu]_{s_2} &&
\KGL \ar[rr]^{\nu} \ar@/_/[uu]_{s} &&
\GVB
\\
E,F & &
f^*E, f^*F & &
E_n, F_n & &
s^*\E'
}
$$
and that the bundles 
${\Et}':=(\id_{\Ct}\times f)^*\Et$ 
and $\G$ are related by the equation
$$
\Et'=(h_*\G(-s_1-s_2))(p_1+p_2)
\quad.
\eqno(1)
$$

It is clear that we have $\nu^*\Theta=\det(R\pi'_*\E')$. 
From the canonical exact sequence of $\Oo_{\Cc'}$-modules
$$
0\to \E' \to g_*\G \to s_*s^*\E' \to 0
$$
and the fact that $R^1g_*\G=0$ (cf. \cite{Knudsen}, Cor 1.5)
we get the canonical isomorphism 
$$
\det(R\pi'_*\E') = \det(R\rho_*\G)\tensor\det(s^*\E')
\quad.
\eqno(2)
$$
The canonical exact sequence of $\Oo_{\Bc}$-modules 
$$
0\to \G(-s_1-s_2)\to\G\to (s_1)_*s_1^*\G\oplus(s_2)_*s_2^*\G \to 0
$$
yields the canonical isomorphism
$$
\det(R\rho_*\G)= \det(R\rho_*\G(-s_1-s_2))\tensor
                (\det s_1^*\G)^{-1}\tensor(\det s_2^*\G)^{-1}
\quad.
\eqno(3)
$$
From equation $(1)$
it follows that there is a canonical exact sequence
of $\Oo_{\Ct\times \KGL}$-modules
$$
0\to h_*\G(-s_1-s_2) \to \Et' \to (p_1)_*p_1^*\Et'\oplus (p_2)_*p_2^*\Et' \to 0
\quad.
$$
From this sequence, the fact that 
$R^1h_*\G(-s_1-s_2)=0$ (cf. \cite{Knudsen}, Cor 1.5)
and the equalities $p_1^*\Et'=f^*E$ and $p_2^*\Et'=f^*F$
it follows that we have canonically:
$$
\det(R\rho_*\G(-s_1-s_2))= \det(R\pit'_*\Et')\tensor 
                          f^*\det E\tensor f^*\det F
\quad.
\eqno(4)
$$
Putting together the identifications $(2)-(4)$ and making use of the
fact that 
$\det(R\pit'_*\Et')=f^*\Thetat$ 
and that
the $\Oo_{\KGL}$-modules 
$s^*\E'$, $s_1^*\G$, $s_2^*\G$, $E_n$, $F_n$ 
are all canonically isomorphic, we finally get
$$
\nu^*\Theta=\det(R\pi'_*\E')=f^*\Thetat\tensor(\det E_n)^{-1}
                             \tensor f^*(\det E)\tensor f^*(\det F)
\quad.
$$
The proposition now follows from \ref{detE_n}.
\end{proof}

Now we will apply the results from \S \ref{cohomology}.
Notation is as in \S \ref{review kgl} and \S \ref{cohomology} with
$S$ replaced by the stack $\VB$ and vector bundles 
$E=p_1^*\Et$, 
$F=p_2^*\Et$, 
on $\VB$ as above.
Furthermore I will write $\PB$ instead of $\Fl$ for the product
$
\text{Fl}(E)\times_{\VB}\text{Fl}(F)
$.
The letters $\PB$ stand of course for {\em parabolic bundles}.
Let $f_{\PB}:\PB\to\VB$ be the canonical projection and
let 
\begin{eqnarray*}
\nu_{I,J} &:=& \nu\comp i_{I,J}: \Ob_{I,J}\to\GVB
\\
f_{I,J} &:=& f\comp i_{I,J}: \Ob_{I,J}\to \VB
\end{eqnarray*}

\begin{proposition}
\label{H^0 of O_IJ}
Let $I,J\subseteq[0,n-1]$ with $\min(I)+\min(J)\geq n$.
Then the following holds:

1.
We have a canonical isomorphism
$$
H^0(\Ob_{I,J}, \nu_{I,J}^*(\Tk))=
\Oplus_{(a,b)\in A_{I,J}(\Dk)}
H^0(\PB,\Tk_{\PB}(a,b))
\quad,
$$
where $\Tk_{\PB}(a,b):=f_{\PB}^*(\Thetat^\kappa)\tensor\Oo(a,b)$
and 
where 
$
A_{I,J}(\Dk)
$
is the set of all $(a,b)\in\Z\times\Z$ with the property that
$b_i=\kappa-a_{n-i+1}$ for $i\in[1,n]$ and
$$
0=a_1=\dots=a_{n-j_1}\leq 
a_{n-j_1+1}\leq\dots\leq a_{i_1}\leq
a_{i_1+1}=\dots=a_n=\kappa
\quad,
$$
where $i_1:=\min(I)$ and $j_1:=\min(J)$. 

2. 
We have $A_{I,J}(\Dk)\subseteq A(\Dk)$
and the following diagram commutes:
$$
\xymatrix{
\text{$H^0(\KGL, \nu^*(\Tk))$}
\ar@{=}[d] 
\ar[r]^{i_{I,J}^{\bullet}} 
&
\text{$H^0(\Ob_{I,J}, \nu_{I,J}^*(\Tk))$}
\ar@{=}[d] 
\\
\text{$
\underset{(a,b)\in A(\Dk)}{\Oplus}
H^0(\PB,\Tk_{\PB}(a,b))
$} 
\ar@<1ex>[r]
&
\text{$
\underset{(a,b)\in A_{I,J}(\Dk)}{\Oplus}
H^0(\PB,\Tk_{\PB}(a,b))
$} 
}
$$
Here, $i_{I,J}^{\bullet}$ is the restriction
morphism induced by the inclusion 
$
i_{I,J}:\Ob_{I,J}\injto\KGL
$
and the lower horizontal arrow is the projection map induced by
the inclusion $A_{I,J}(\Dk)\subseteq A(\Dk)$.
\end{proposition}

\begin{proof}
By Theorem \ref{bundle decomposition} we have 
a canonical decomposition
$$
(f_{I,J})_*\Dk=\Oplus_{(a,b)\in A_{I,J}(\Dk)}(f_{\PB})_*\Oo(a,b)
\quad.
$$ 
This, together with \ref{pullback of Theta}
implies the isomorphism stated in the first part of the proposition.
The concrete description of the set $A_{I,J}(\Dk)$ follows easily
from definition \ref{A}.
The second part of the proposition
is immediate from \ref{bundle decomposition}.2.
\end{proof}

The main result of this paper is the following

\begin{theorem}
\label{main}
There is a canonical isomorphism
$$
H^0(\GVB,\Tk)\isomto
\Oplus_{(a,b)\in A'}
H^0(\PB,\Tk_{\PB}(a,b))
\quad,
$$
where $A'$ is the set of all $(a,b)\in\Z^n\times\Z^n$ 
with $0\leq a_1\leq\dots\leq a_n\leq\kappa-1$ and $b_i=\kappa-a_{n-i+1}$
for $i\in[1,n]$.
\end{theorem}

The remaining of this section is devoted to the proof of theorem
\ref{main}.

\begin{definition}
For $p,q\in[0,n]$ with $p+q\geq n$ we set
\begin{eqnarray*}
A_{p,q} &:=& \{(a,b)\in A(\Dk)\ |\ 
             \text{$a_i=0$ for $i\in[1,n-q]$ and
                   $a_i=\kappa$ for $i\in[p+1,n]$} \} \\
A'_{p,q} &:=& \{(a,b)\in A_{p,q}\ |\ 
             \text{$a_i\leq\kappa-1$ for $i\in[1,p]$} \} \\
V_{p,q} &:=& \Oplus_{(a,b)\in A_{p,q}}H^0(\PB,\Tk_{\PB}(a,b)) \\
V'_{p,q} &:=& \Oplus_{(a,b)\in A'_{p,q}}H^0(\PB,\Tk_{\PB}(a,b)) 
\end{eqnarray*}
\end{definition}

\begin{remark}
\label{remark}
1.
Let $I,J\subseteq[0,n-1]$, let $p:=\min(I)$ and $q:=\min(J)$, 
and assume $p+q\geq n$.
Then by \ref{H^0 of O_IJ} we have  a canonical 
isomorphism
$$
V_{p,q}=H^0(\Ob_{I,J},\nu_{I,J}^*\Tk)
\quad.
$$

2.
For $p,p',q,q'\in[0,n]$ with $p\leq p'$, $q\leq q'$ and $p+q\geq n$
we have $A_{p,q}\subseteq A_{p',q'}$.
Furthermore we have $A_{n,n}=A(\Dk)$, and 
$A'_{n,n}=A'$ is the set which appears in theorem \ref{main}.

3.
Let $p,q\in[0,n]$ with $p+q\geq n$.
Then $A_{p,q}$ is the disjoint union of the sets $A'_{i,q}$, where
$i$ runs through $[n-q,p]$. Therefore we have 
$
V_{p,q}=\Oplus_{i=n-q}^pV'_{i,q}
$.
It follows that $V_{p,q}=V'_{p,q}$, if $p+q=n$ and 
$V_{p,q}=V'_{p,q}\oplus V_{p-1,q}$, if $p+q>n$.
\end{remark}

\begin{definition}
\label{morphisms}
1.
Let $p,p',q,q'\in[0,n]$ with $p\leq p'$, $q\leq q'$ and $p+q\geq n$.
Then we denote by 
$$
\sigma^{p',q'}_{p,q}: V_{p',q'}\To V_{p,q}
\qquad\text{and}\qquad
\tau^{p',q'}_{p,q}: V'_{p',q'}\To V'_{p,q}
$$
the projection morphisms induced by the inclusions
$
A_{p,q}\injto A_{p',q'}
$
and
$
A'_{p,q}\injto A'_{p',q'}
$
respectively.

2.
For $p,q\in[0,n]$ with $p+q\geq n$ we denote by 
$$
\pi_{p,q}:V_{p,q}\To V'_{p,q}
$$
the projection induced by the inclusion $A'_{p,q}\injto A_{p,q}$.

3.
Let $p\in [0,n-1]$.  
We denote by
$$
\beta^{p,n}_{n,p}:
V_{p,n}
=
H^0(Z_p,\nu_{Z_p}^*\Theta^{\kappa})
\Isomto
H^0(Y_p,\nu_{Y_p}^*\Theta^{\kappa})
=
V_{n,p}
$$
the isomorphism 
induced on cohomology by 
$
\beta_p:
Y_p\isomto Z_p
$
(cf. \ref{H^0 as subspace})
via the identification \ref{remark}.1.
For convenience, we  define $\beta^{n,n}_{n,n}$ to be the identity
morphism on $V_{n,n}$.
\end{definition}

\begin{remark}
\label{remark2}
1.
By \ref{H^0 of O_IJ}.2 the morphisms 
$i_{Z_p}^\bullet$,
$i_{Y_p}^\bullet$ and
$\beta_p^\bullet$
from \ref{H^0 as subspace} 
are equal (via the identification \ref{remark}.1)
to the morphisms
$\sigma^{n,n}_{p,n}$,
$\sigma^{n,n}_{n,p}$ and
$\beta^{p,n}_{n,p}$ 
respectively. Thus by \ref{H^0 as subspace} the space 
$H^0(\GVB,\Tk)$ can be identified with the subspace
of all $\theta\in V_{n,n}$ which have the property that
$$
\beta^{p,n}_{n,p}\ \sigma^{n,n}_{p,n}\ \theta\ =\ \sigma^{n,n}_{n,p}\ \theta
$$
for every $p\in[0,n-1]$.

2.
The following equalities are trivially verified:
$$
\sigma^{p',q'}_{p,q}\comp\sigma^{p'',q''}_{p',q'}=\sigma^{p'',q''}_{p,q}
\qquad\text{and}\qquad
\pi_{p,q}\comp\sigma^{p',q'}_{p,q}=\tau^{p',q'}_{p,q}\comp\pi_{p',q'}
\quad.
$$
\end{remark}

\begin{lemma}
\label{V' to V'}
Let $p\in[0,n]$.
Then the isomorphism
$$
\beta^{p,n}_{n,p}: V_{p,n}\isomto V_{n,p}
$$
maps the subspace $V'_{p,n}$ onto the subspace $V'_{n,p}$.
\end{lemma}

\begin{proof}
For $p=n$ the assertion is trivial, so let $p\in[0,n-1]$.
Assume for a moment that there exist line bundles 
\begin{eqnarray*}
\M=(\Tensor_{i=0}^{n-1}M_i^{m_i})
\qquad\text{and}\qquad
\M'=(\Tensor_{i=0}^{n-1}M_i^{m'_i})
\end{eqnarray*}
on $\KGL$ with the following properties:
\begin{enumerate}
\item
The $m_i$ and $m'_i$ are non-negative and
$m_i=0$ for $i\in[p,n-1]$.
\item
$
\beta_p^*
i_{Z_p}^*
(\M,\mu_0^{m_0}\tensor\dots\tensor\mu_{n-1}^{m_{n-1}})
=
i_{Y_p}^*
(\M',\mu_0^{m'_0}\tensor\dots\tensor\mu_{n-1}^{m'_{n-1}})
$
\item
\begin{eqnarray*}
A'_{p,n} &=&
A_{\{p\},\emptyset}(\Dk\tensor\M^{-1})
\\
A'_{n,p} &=&
A_{\emptyset,\{p\}}(\Dk\tensor(\M')^{-1})
\end{eqnarray*}
\end{enumerate}
By property 1 multiplication with 
$\mu_0^{m_0}\tensor\dots\tensor\mu_{n-1}^{m_{n-1}}$
induces an injection
$$
\nu_{Z_p}^*\Tk\tensor i_{Z_p}^*\M^{-1} 
\injto 
\nu_{Z_p}^*\Tk
\quad,
$$
It follows from property 2 that application of the functor 
$\beta_p^*$
to this injection yields the injection
$$
\nu_{Y_p}^*\Tk\tensor i_{Y_p}^*(\M')^{-1} 
\injto 
\nu_{Y_p}^*\Tk
$$
induced by multiplication with 
$\mu_0^{m'_0}\tensor\dots\tensor\mu_{n-1}^{m'_{n-1}}$.
Therefore we have a commutative diagram:
$$
\xymatrix{
H^0(Z_p,\nu_{Z_p}^*\Tk)
\ar[r]^{\isomorph}
&
H^0(Y_p,\nu_{Y_p}^*\Tk)
\\
H^0(Z_p,
    \nu_{Z_p}^*\Tk\tensor i_{Z_p}^*\M^{-1})
\ar@{^(->}[u] \ar[r]^{\isomorph}
&
H^0(Y_p,
    \nu_{Y_p}^*\Tk\tensor i_{Y_p}^*(\M')^{-1}) 
\ar@{^(->}[u]
}
$$
where the horizontal arrows are induced by the isomorphism 
$\beta_p$.
By property 3 and \ref{bundle decomposition} 
this diagram may be identified with a diagram of the form
$$
\xymatrix@C=8ex@R=2ex{
V_{p,n} \ar[r]^{\beta^{p,n}_{n,p}}
&
V_{n,p}
\\
V'_{p,n} \ar[r]^{\isomorph} \ar@{^(->}[u]
&
V'_{n,p} \ar@{^(->}[u]
}
$$
where the vertical arrows are induced by the inclusions
$A'_{p,n}\injto A_{p,n}$ 
and
$A'_{n,p}\injto A_{n,p}$.
This clearly implies the lemma.

Thus it remains only to prove the existence of $\M$ and $\M'$.
For this, let $c,c'\in\Z^n$ be defined by
$$
c_i:=
\left\{
\begin{array}{ll}
\kappa-1  & \text{if $i\in[1,p]$} \\
\kappa  & \text{if $i\in[p+1,n]$} 
\end{array}
\right.
\qquad\text{and}\qquad
c'_i:=
\left\{
\begin{array}{ll}
0  & \text{if $i\in[1,n-p]$} \\
\kappa-1  & \text{if $i\in[n-p+1,n]$} 
\end{array}
\right.
$$
and let
$$
m_i:=\kappa(n-i)-\sum_{j=i+1}^nc_j
\qquad\text{and}\qquad
m'_i:=\kappa(n-i)-\sum_{j=i+1}^nc'_j
$$
for $i\in[0,n]$.
A simple calculation shows that the line bundles $\M$ and $\M'$
formed with this choice of $m_i$ and $m'_i$ have the properties 
1 and 3.
Property 2 follows easily from proposition \ref{isomorphism}.3.
\end{proof}

After these preparations we now come to the proof of theorem \ref{main}.
I claim that the composite morphism
$$
\xymatrix{
\text{$H^0(\GVB,\Tk)$} \ar@{^(->}[r] &
\text{$H^0(\KGL,\nu^*\Tk) = V_{n,n}$} \ar[r]^(.4){\pi_{n,n}} &
\text{$V'_{n,n}=\Oplus_{(a,b)\in A'}H^0(\PB,\Tk_{\PB}(a,b))$}
}
$$
is an isomorphism.

To prove injectivity, let $\theta\in V_{n,n}$ be an element in the
kernel of $\pi_{n,n}$, which satisfies the 
condition stated in \ref{remark2}.1. We have to show that $\theta=0$.

Since 
$$
V_{n,n}=\Oplus_{p=0}^n V'_{p,n}
$$
it suffices to show that 
$\pi_{p,n}\ \sigma^{n,n}_{p,n}\theta=0$ for all $p\in[0,n]$.
We do this by induction on $p$.

For $p=0$ we have
$
\beta^{0,n}_{n,0}\ \sigma^{n,n}_{0,n}\ \theta =
\sigma^{n,n}_{n,0}\ \theta=
\pi_{n,0}\ \sigma^{n,n}_{n,0}\ \theta=
\tau^{n,n}_{n,0}\ \pi_{n,n}\ \theta=0
\quad,
$
which implies 
$
\pi_{0,n}\ \sigma^{n,n}_{0,n}\ \theta = 
\sigma^{n,n}_{0,n}\ \theta = 0
$.
Now let $p>0$ and assume 
$\pi_{q,n}\ \sigma^{n,n}_{q,n}\ \theta=0$ for all $q\in[0,p-1]$.
This implies that $\sigma^{n,n}_{p,n}\ \theta$ is in fact contained in
$V'_{p,n}$. Therefore by \ref{V' to V'} we have that 
$
\beta^{p,n}_{n,p}\ \sigma^{n,n}_{p,n}\ \theta =
\sigma^{n,n}_{n,p}\ \theta 
$
is contained in $V'_{n,p}$.
But this implies
$
\beta^{p,n}_{n,p}\ \sigma^{n,n}_{p,n}\ \theta = 
\pi_{n,p}\ \sigma^{n,n}_{n,p}\ \theta =
\tau^{n,n}_{n,p}\ \pi_{n,n}\ \theta = 0
$
and thus 
$
\pi_{p,n}\ \sigma^{n,n}_{p,n}\ \theta=
\sigma^{n,n}_{p,n}\ \theta=0
$.

It remains to prove surjectivity.
Let $\theta'$ be an element of $V'_{n,n}$.
For $p\in[0,n]$ let $\theta'_p\in V'_{p,n}$
be defined inductively by the property
$$
\beta^{p,n}_{n,p}\ \theta'_p =
\tau^{n,n}_{n,p}\ \theta' -
\sum_{q=0}^{p-1}\ \pi_{n,p}\ \beta^{p,n}_{n,p}\ \theta'_q
$$
and let $\theta_p:=\sum_{q=0}^{p}\ \theta'_q$ and $\theta:=\theta_n$.
Clearly, we have $\pi_{n,n}\ \theta=\theta'_n=\theta'$. Therefore it
suffices to show that $\theta$ is an element of $H^0(\GVB,\Tk)$.

By \ref{remark2}.1 this amounts to proving that 
$
\beta^{q,n}_{n,q}\ \theta_q=
\sigma^{n,n}_{n,q}\ \theta
$
for $q\in[0,n]$.
Since we have
$$
V_{n,q}=\Oplus_{p=n-q}^nV'_{p,q}
\quad,
$$
this is equivalent to the statement that
$$
\pi_{p,q}\ \sigma^{n,q}_{p,q}
\beta^{q,n}_{n,q}\ \theta_q=
\pi_{p,q}\ \sigma^{n,n}_{p,q}\ \theta
$$
for all $p,q\in[0,n]$ with $p+q\geq n$.

But this equality is clear, since
\begin{eqnarray*}
\pi_{p,q}\ \sigma^{n,q}_{p,q}
\beta^{q,n}_{n,q}\ \theta_q 
&=&
\tau^{n,q}_{p,q}\ \pi_{n,q}
\beta^{q,n}_{n,q}(\theta'_q+\theta_{q-1})= 
\\
&=&
\tau^{n,q}_{p,q}\ \pi_{n,q}
(\tau^{n,n}_{n,q}\ \theta'-
 \pi_{n,q}\ \beta^{q,n}_{n,q}\ \theta_{q-1}+
 \beta^{q,n}_{n,q}\ \theta_{q-1})=
\\
&=&
\tau^{n,n}_{p,q}\ \theta'=
\pi_{p,q}\ \sigma^{n,n}_{p,q}\ \theta
\end{eqnarray*}

\begin{remark}
The decomposition given in Theorem \ref{main} is not symmetric
with respect to the two points $p_1$ and $p_2$.
Indeed, interchanging the role of the two points means interchanging
$a$ and $b$ in $(a,b)$, but then the set $A'$ changes
to the set $A''$ consisting of all $(b,a)\in\Z\times\Z$ with
$1\leq b_1\leq\dots\leq b_n\leq\kappa$ and $a_i=\kappa-b_{n-i+1}$.
At first sight this seems strange, since the decomposition should
certainly not depend on how we numerate the points $p_1$ and $p_2$.
The answer to this riddle is that in our proof of Theorem \ref{main}
we have made a choice between two possibilities.
In fact, one can equally well show that the composite
morphism
$$
\xymatrix{
H^0(\GVB,\Tk) \ar@{^(->}[r] &
\text{$
\Oplus_{(a,b)\in A(\Dk)}
H^0(\PB,\Tk_{\PB}(a,b))
$}
\ar[r] &
\text{$
\Oplus_{(a,b)\in A''}
H^0(\PB,\Tk_{\PB}(a,b))
$}
}
$$
is an isomorphism, where the last arrow is the projection morphism
induced by the inclusion $A''\subset A(\Dk)$.
\end{remark}

\begin{remark}
\label{chi}
Let $\X$ be an algebraic $k$-stack and let $\F$ be a sheaf on the
smooth-\'etale site of $\X$ (cf. \cite{LM} \S 12). 
By definition (loc. cit. (12.5.3)), 
the set of global sections of $\F$ is the set
of all families $s_{(U,u)}$ of sections of $\F$ over 
$(U,u)\in\ob\ \Liset(\X)$ such that 
$
\res_{\varphi}s_{(V,v)} = s_{(U,u)}
$
for all arrows $\varphi:(U,u)\to(V,v)$ in $\Liset(\X)$.
Now assume in particular that $\F$ is an $\Oo_{\X}$-module
and that for each object $(U,u)\in\ob\ \Liset(\X)$ there
is a homomorphism $k^{\times}\to\Aut_{\Liset(\X)}(U,u)$, $a\mapsto\varphi_a$.
Assume furthermore that there is a number $\chi\in\Z$ 
such that for each $(U,u)$ and $a\in k^{\times}$ the morphism
$\res_{\varphi_a}:\F(U,u)\to\F(U,u)$
is multiplication with the $\chi$-th power of $a$.
Then it is clear that unless $\chi=0$, the only global
section of $\F$ is the zero section.

The stack $\GVB$ is the disjoint union of open closed substacks
$\GVB_d$ parametrizing Gieseker vector bundles of degree $d$ ($d\in\Z$).
By the above consideration it follows that $H^0(\GVB_d,\Theta^{\kappa})$ 
vanishes unless the Euler characteristic
$$
\chi=d+n(1-g)
$$
of a bundle of rank $n$ and degree $d$ on a curve of genus $g$
is zero, i.e. unless $d=n(g-1)$.
Therefore we have 
$H^0(\GVB,\Theta^{\kappa})=H^0(\GVB_{n(g-1)},\Theta^{\kappa})$.
A similar remark applies to the groups 
$
H^0(\PB,\Theta^{\kappa}_{\PB}(a,b))
$
appearing on the right hand side of the isomorphism in
\ref{main}.
\end{remark}

%-----------------------------------------------------------------------------
\section{Degeneration}

Let $B$ be the spectrum of a discrete valuation ring and
let $C\to B$ be a projective relative curve of genus 
$g\geq 1$ over $B$, whose generic fiber $C_{\eta}$
is smooth and whose special fiber $C_0$ is irreducible with one ordinary
double point $p$. 
In \cite{degeneration} I have shown that there is a flat algebraic 
moduli stack $\GVB(C/B)$ over $B$ whose generic fiber $\VB(C_{\eta})$
parametrizes
vector bundles on $C_{\eta}$ and whose special fiber $\GVB(C_0)$
parametrizes 
Gieseker vector bundles of rank $n$ on $C_0$.
Let 
$$
(\pi_{C/B}:\Cc_{C/B}\to\GVB(C/B)\ ,\ \E_{C/B})
$$ 
be the universal Gieseker vector bundle over $\GVB(C/B)$
and let 
$$
\Theta(C/B):=\det R(\pi_{C/B})_*\E_{C/B}
$$
be the determinant line bundle on $\GVB(C/B)$.

The rest of this section is dedicated to the proof
of the following result which shows that the model
$\GVB(C/B)$ of $\VB(C_{\eta})$ defines the 
``correct selection rules'' for generalized
theta functions.

\begin{theorem}
The $B$-module
$H^0(\GVB(C/B),\Theta(C/B)^{\kappa})$ is locally free of finite rank.
\end{theorem}

It is clear that there is a decomposition into a disjoint union:
$$
\GVB(C/B)=\bigsqcup_{d\in\Z}\GVB_d(C/B)
\quad,
$$
where $\GVB_d(C/B)$ is the open substack of $\GVB(C/B)$
parametrizing vector bundles of degree $d$.
As in \ref{chi} it follows that
there are no non-vanishing sections of $\Theta(C/B)^{\kappa}$
over $\GVB_d(C/B)$ unless $d=n(g-1)$.
From now on we let $d:=n(g-1)$ and we will 
restrict our attention to the open substack
$\GVB_d(C/B)$. Its closed and special fiber over $B$ will
be denoted by $\GVB_d(C_0)$ and $\VB_d(C_{\eta})$ and the restriction
of $\Theta(C/B)$ to these by $\Theta(C_0)$ and $\Theta(C_{\eta})$
respectively.

We have to show that for $\kappa\geq 1$ the vector spaces
$H^0(\GVB_d(C_0),\Theta(C_0)^{\kappa})$ and  
\linebreak
$H^0(\VB_d(C_{\eta}),\Theta(C_{\eta})^{\kappa})$ 
are finite dimensional and that the equality
$$
\dim H^0(\GVB_d(C_0),\Theta(C_0)^{\kappa})
=
\dim H^0(\VB_d(C_{\eta}),\Theta(C_{\eta})^{\kappa})
\eqno(1)
$$
holds.
Since the stacks $\GVB_d(C_0)$ and $\VB_d(C_{\eta})$
are not separated and not of finite type over their 
respective base fields, we cannot apply general results
from \cite{LM} \S 15 or \cite{F3} to prove finite dimensionality
of cohomology. By the same reason we cannot argue by cohomological
flatness to show that the dimensions coincide,
even if we could assume that the relevant theorems in
\cite{EGA} III hold in the context of Artin stacks.
Our strategy therefore is to compute the dimensions
of
$H^0(\GVB_d(C_0),\Theta(C_0)^{\kappa})$
and of
$H^0(\VB_d(C_{\eta}),\Theta(C_{\eta})^{\kappa})$ 
individually by relating them to the dimensions of spaces of 
generalized theta functions for (parabolic) $\SL_n$-bundles and
using the Verlinde formula.

\begin{definition}
\def\theenumi{\roman{enumi}}
\def\labelenumi{(\theenumi)}
\begin{enumerate}
\item
We fix once and for all a line bundle $L_{\eta}$ of degree $d$ on 
$C_{\eta}$.
We denote by $\SVB(C_{\eta})$ 
the closed substack of $\VB_d(C_{\eta})$ 
which parametrizes vector bundles whose determinant is isomorphic
to $L_{\eta}$. We write $\Theta_{\SVB}(C_{\eta})$ for  the restriction
of $\Theta(C_{\eta})$ to $\SVB(C_{\eta})$.
\item
Also we fix once and for all a line bundle $\tilde{L}_0$ 
of degree $d$ on $\Ct$.
As in the previous paragraphs we let $\VB$ denote the moduli stack
of rank $n$ vector bundles on $\Ct$.
Let $\VB_d=\VB_d(\Ct)$ be the open substack of $\VB$ which parametrizes
bundles of degree $d$ and let $\SVB=\SVB(\Ct)$ be the closed substack
of $\VB$ which parametrizes vector bundles whose determinant is isomorphic
to $\tilde{L}_0$.
Recall from \S \ref{section decomposition} that
$\PB$ denotes the stack of vector bundles on $\Ct$ together
with full flags in the fibers over the points $p_1$ and $p_2$. 
We define
$$
\PB_d:=\PB\times_{\VB}\VB_d
\quad,\quad
\SPB:=\PB\times_{\VB}\SVB
$$
and for $(a,b)\in\Z^n\times\Z^n$ and $\kappa\in \Z$
we denote by 
$\Theta_{\SPB}$, $\Oo_{\SPB}(a,b)$ and
$\Theta^{\kappa}_{\SPB}(a,b)$ 
the restriction of the line bundles
$\Theta_{\PB}=f_{\PB}^*\Thetat$, $\Oo_{\PB}(a,b)$ and 
$\Theta^{\kappa}_{\PB}(a,b)=f_{\PB}^*\Thetat^{\kappa}\tensor\Oo_{\PB}(a,b)$ 
respectively
to the closed substack $\SPB$ of $\PB$. 

\end{enumerate}
\end{definition}

\begin{proposition}
\label{donagi}
\def\theenumi{\roman{enumi}}
\def\labelenumi{(\theenumi)}
\begin{enumerate}
\item
The dimensions of the vector spaces 
$H^0(\SVB(C_{\eta}),\Theta_{\SVB}(C_{\eta})^{\kappa})$
and
%\linebreak
$H^0(\VB_d(C_{\eta}),\Theta(C_{\eta})^{\kappa})$
are finite and 
we have
$$
\dim H^0(\VB_d(C_{\eta}),\Theta(C_{\eta})^{\kappa})
=
\left(\frac{\kappa}{n}\right)^g\cdot 
\dim H^0(\SVB(C_{\eta}),\Theta_{\SVB}(C_{\eta})^{\kappa})
\quad.
$$ 
\item
Let $(a,b)\in\Z^n\times\Z^n$
and let $(a',b')\in\Z^n\times\Z^n$ be defined by
$a'_i:=a_i-a_1$, $b'_i:=b_i-b_1$.
Then the vector spaces 
$H^0(\PB_d,\Theta^\kappa(a,b))$
and 
$H^0(\SPB,\Theta^{\kappa}_{\SPB}(a',b'))$
are finite dimensional and we have
$$
\dim H^0(\PB_d,\Theta^\kappa(a,b))=
\left(\frac{\kappa}{n}\right)^{g-1}\cdot 
\dim H^0(\SPB,\Theta^{\kappa}_{\SPB}(a',b'))
\quad.
$$
\end{enumerate}
\end{proposition}

\begin{proof}
The finiteness of  
$H^0(\SVB(C_{\eta}),\Theta_{\SVB}(C_{\eta})^{\kappa})$
and
$H^0(\SPB,\Theta^{\kappa}_{\SPB}(a',b'))$
follows from the interpretation of these vector spaces
as spaces of conformal blocks (cf. \cite{BL} and \cite{Pauly}).

The analogous equality to $(i)$ in the context of coarse moduli
spaces of semi-stable bundles has been proved by Donagi and Tu
in 
\cite{DT}. 
Since their proof works almost identically 
in our situation, we will concentrate on the second equation,
of which $(i)$ can anyway be considered a special case (namely
the case $a=b=0$).

For the proof of $(ii)$ consider the following Cartesian diagram
$$
\xymatrix@C=16ex{
\SPB\times\J_0 \ar[r]^{\sigma} \ar[d]_{f_{\SPB}\times\id} 
&
\PB_d \ar[d]^{f_{\PB}} 
\\
\SVB\times\J_0 \ar[r]^{\tau}_{E,L\mapsto E\tensor L} \ar[d]_{\pr}
&
\VB_d \ar[d]^{\det}
\\
\J_0 \ar[r]^{\rho}_{L\mapsto L^n\tensor \tilde{L}_0}
&
\J_d
}
$$
where $\J_0$ and $\J_d$ are the moduli stacks of line bundles on
$\Ct$ of degree $0$ and $d$ respectively. 
Similarly as in \cite{DT}, the sought-for equality follows from
computing 
the space of global sections of the line bundle
$\sigma^*\Theta^{\kappa}_{\PB}(a,b)$ on $\SPB\times\J_0$
in two different ways.

{\bf First way:}
Let $\E_{\PB}$, $\E_{\J_0}$, etc. 
be the universal vector bundle on $\Ct\times\PB$, $\Ct\times\J_0$, etc..
Let $\pi_{\PB}$, $\pi_{\J_0}$, etc.
be the projection from $\Ct\times\PB$, $\Ct\times\J_0$ etc.
to the second factor
and denote by
$p_1,p_2$ be the sections of $\pi_{\PB}$, $\pi_{\J_0}$ etc. 
induced by the points
$p_1,p_2\in \Ct$.
We have
\begin{eqnarray*}
\Theta^{\kappa}_{\PB}(a,b) 
&=&
(\det p_2^*\E_{\PB} \tensor \det R\pi_{\PB *}\E_{\PB})^{\kappa}
\\
&\tensor & 
(\det p_2^* \E_{\PB})^{a_1+b_1-\kappa}\tensor\Oo_{\PB}(a',b')
\\
& \tensor& 
(\det p_1^*\E_{\PB}\tensor(\det p_2^*\E_{\PB})^{-1})^{a_1}
\end{eqnarray*}

Since 
$
\det p_2^*\E_{\PB} \tensor \det R\pi_{\PB *}\E_{\PB} =
f_{\PB}^* \det R \pit_* \Et(-p_2)
$,
it follows as in \cite{DT}, Cor 6 that we have
$$
\sigma^*(\det p_2^*\E_{\PB} \tensor \det R\pi_{\PB *}\E_{\PB}) =
\Theta_{\SVB} \boxtimes \Theta_N^n
\quad,
$$
where $N\in\Pic^{g-2}(\Ct)$ is an $n$-th root of 
$\tilde{L}_0\tensor\Oo_{\Ct}(-np_2)$ and
where $\Theta_N$ is the line bundle 
$\det R\pi_{\J_0 *} (\E_{\J_0}\tensor N)$
on $\J_0$.
It follows directly from the definitions that we have
$$
\sigma^*
((\det p_2^* \E_{\PB})^{a_1+b_1-\kappa}\tensor\Oo_{\PB}(a',b'))
=
\Oo_{\SPB}(a',b') \boxtimes 
(p_1^*\E_{\J_0}\tensor p_2^*\E_{\J_0}^{-1})^{\sum a'_i}
$$
and
$$
\sigma^*
(\det p_1^*\E_{\PB}\tensor(\det p_2^*\E_{\PB})^{-1})
=
\Oo_{\SPB}\boxtimes(p_1^*\E_{\J_0}\tensor p_2^*\E_{\J_0}^{-1})^n
\quad.
$$

Summarizing, we have 
$$
\sigma^*\Theta_{\PB}^{\kappa}(a,b)=
\Theta_{\SPB}^{\kappa}(a',b')
\boxtimes
\left(
\Theta_N^{\kappa n}\tensor
(p_1^*\E_{\J_0}\tensor p_2^*\E_{J_0}^{-1})^{\sum a_i}
\right)
\quad.
$$
Using the techniques from \cite{DT} \S5 we see
that
$$
\Theta_N^{\kappa n}\tensor
(p_1^*\E_{\J_0}\tensor p_2^*\E_{J_0}^{-1})^{\sum a_i}
=
\tau_{M}^*\Theta_N^{\kappa n}
\quad,
$$
where $M\in \Pic^0(\Ct)$ is a $\kappa n$-th root of 
the line bundle $\Oo_{\Ct}((\sum_{i=1}^n a_i)(p_2-p_1))$ and
$\tau_M:\J_0\to\J_0$ is the translation by $M$.
Thus we have
$$
\sigma^*\Theta_{\PB}^{\kappa}(a,b)=
\Theta_{\SPB}^{\kappa}(a',b')
\boxtimes
\tau_{M}^*\Theta_N^{\kappa n}
$$
and consequently
$$
H^0(\SPB\times\J_0,\sigma^*\Theta_{\PB}^{\kappa}(a,b))
\isomorph
H^0(\SPB,\Theta_{\SPB}^{\kappa}(a',b'))
\tensor
H^0(\J_0,\Theta_N^{\kappa n})
\quad.
\eqno(*)
$$

{\bf Second way:}
The morphism $\sigma$ is a Galois covering with 
Galois group $G$ the subgroup of $n$-torsion points
of $\Pic^0(\Ct)$.
As in \cite{DT}, Prop. 4. and Lemma 7. we have
$$
\sigma_*\Oo_{\SPB\times\J_0}=\bigoplus_{\lambda\in\hat{G}}L_{\lambda}
\quad,
$$
where $\hat{G}$ is the character group of $G$ and for 
each $\lambda\in\hat{G}$ we have
$$
L_{\lambda}=(\det\comp f_{\PB})^*N_{\lambda}
$$
for some line bundle $N_{\lambda}$ of degree zero on $\J_d$.
By the projection formula we have
$$
\sigma_*\sigma^*\Theta_{\PB}^{\kappa}(a,b)
=
\bigoplus_{\lambda\in\hat{G}}
\Theta_{\PB}^{\kappa}(a,b)\tensor L_{\lambda}
\quad.
$$
As in \cite{DT} \S5 it follows that 
$$
\Theta_{\PB}^{\kappa}(a,b)\tensor L_{\lambda}
=
\tau_{M_{\lambda}}^*\Theta_{\PB}^{\kappa}(a,b)
\quad,
$$
for some line bundle $M_{\lambda}\in\Pic^0(\Ct)$, 
where $\tau_{M_{\lambda}}:\PB\to\PB$ is the isomorphism
which sends a parabolic bundle $E$ to $E\tensor M_{\lambda}$.
Therefore we have
$$
H^0(\SPB\times\J_0,\sigma^*\Theta_{\PB}^{\kappa}(a,b))
=
H^0(\PB,\sigma_*\sigma^*\Theta_{\PB}^{\kappa}(a,b))
\isomorph
\bigoplus_{\lambda\in\hat{G}}
H^0(\PB,\Theta_{\PB}^{\kappa}(a,b))
\quad.
\eqno(**)
$$

{\bf Conclusion:}
The sought-for equation follows from $(*)$ and $(**)$ together 
with the
fact that the group $\hat{G}$ is of order $n^{2(g-1)}$ and
the fact that we have
$$
\dim H^0(\J_0,\Theta^{\kappa n}_N) = (\kappa n)^{g-1}
\quad.
$$
This last equality is well known for theta functions on the Jacobian
variety $J_0$, but since we are dealing here with the stack $\J_0$, it requires
some further justification.
For this let $K$ be the open subscheme of a Quot-scheme which parametrizes
invertible quotients $L$ of $\Oo_{\Ct}((2-2g)p_1)^{g-1}$ such that
$H^1(\Ct,L((2g-2)p_1))=0$ and such that the induced morphism
$k^{g-1}\to H^0(\Ct,L((2g-2)p_1))$ is an isomorphism.
Let $\E_K$ be the universal quotient bundle on $\Ct\times K$ and let
$K'$ be the complement of the zero section of (the total space of) the 
line bundle 
$p_1^*\E_K$ on $K$.
Then $\GL_{g-1}$ operates in an obvious way on
$K$ and $K'$ such that the center $\Gm\subset\GL_{g-1}$ operates
trivially on $K$ and 
the projection $K'\to K$ is a $\Gm$-torsor.
We have
$$
\J_0=[K/\GL_{g-1}]
\qquad\text{and}\qquad
J_0=[K'/\GL_{g-1}]
\quad.
$$
Since at each point of $K$ the group $\Gm$ operates trivially on
the fiber of the line bundle $\det R\pi_{K*}(\E_K\tensor N)$
we have
\begin{eqnarray*}
H^0(\J_0,\Theta^{\kappa n}_N)
&=&
H^0(K,(\det R\pi_{K*}(\E_K\tensor N))^{\kappa n})^{\GL_{g-1}}
\\
&=&
H^0(K',(\det R\pi_{K'*}(\E_{K'}\tensor N))^{\kappa n})^{\GL_{g-1}}
=
H^0(J_0,\Theta^{\kappa n}_N)
\quad,
\end{eqnarray*}
where the superscript $\GL_{g-1}$ means taking invariants 
under this group.
\end{proof}

Since by Theorem \ref{main} we have
$$
\dim H^0(\GVB_d(C_0),\Theta(C_0)^{\kappa})
=
\sum_{(a,b)\in A'}
\dim H^0(\PB_d,\Theta^{\kappa}(a,b))
$$
it follows from \ref{donagi} that the vector space
$H^0(\GVB_d(C_0),\Theta(C_0)^{\kappa})$ is finite and that the 
equality (1)
is equivalent to the equality
$$
\sum_{(a',b')\in \SA'}
(\kappa-a'_n)\cdot
\dim H^0(\SPB,\Theta^{\kappa}_{\SPB}(a',b'))
=
\frac{\kappa}{n}
\cdot
\dim H^0(\SVB(C_{\eta}),\Theta_{\SVB}(C_{\eta})^{\kappa})
\quad,
\eqno(2)
$$
where $\SA'$ is the set of all $(a',b')\in\Z^n\times\Z^n$ with
the property that
$
0=a'_1\leq a'_2\leq\dots\leq a'_n\leq\kappa
$
and
$
b'_i=a'_n-a'_{n-i+1}
$.

It is well known that the dimensions of the vector
spaces
$
H^0(\SPB,\Theta^{\kappa}_{\SPB}(a',b'))
$
and
$
H^0(\SVB(C_{\eta}),\Theta_{\SVB}(C_{\eta})^{\kappa})
$
are given by the Verlinde formula.
To write down the formulas explicitly, we need to introduce 
some notation.

Let 
$$
P=\left(\bigoplus_{i=1}^n \Z \epsilon_i\right)/
\left(\sum_{i=1}^n\epsilon_i\right)
$$
be the weight lattice of $\sln$.
Let $(\ \  |\ \ ):P\times P\to Z[1/n]$ be the normalized Killing form
defined by
$$
(\epsilon_i|\epsilon_j):=\delta_{i,j}-\frac{1}{n}
\quad.
$$
Let
$$
R_+:=\{\epsilon_i-\epsilon_j\ |\ 1\leq i<j\leq n\}
$$
be the set of positive roots of $\sln$.
Let
\begin{eqnarray*}
\theta &:=& \epsilon_1 - \epsilon_n \\
\rho &:=& \sum_{i=1}^n(n-i)\epsilon_i
\end{eqnarray*} 
the highest root and the half sum of all positive roots respectively.
Let
$
P_+:=\{\lambda\in P\ |\ 
\text{
$(\alpha|\lambda)\geq 0$
for all $\alpha\in R_+$
} \}
$
be the set of dominant weights 
and let
$$
P_{\kappa}:=\{\lambda\in P_+\ |\ 
(\theta|\lambda)\leq \kappa
 \}
\quad.
$$
Recall that $P_+$ parametrizes the finite dimensional representations
of $\sln$.
The Weyl group $W=S_n$ operates on $P$ by permuting the 
generators $\epsilon_i$. 
Let $w_0:j\mapsto n-j+1$ be the longest element in $W$.
Then $\lambda\mapsto\lambda^*:=-w_0\lambda$ is an involution
of the set $P_+$ (and $P_{\kappa}$), which corresponds to
taking the dual representation.

For $\lambda,\mu\in P$ we define
the complex number
$$
J(\lambda,\mu):=
\sum_{w\in W}
\sign(w)
\exp\left(
\frac{2\pi i}{n+\kappa}
(w(\lambda)|\mu)
\right)
\quad.
$$

\begin{proposition}
\label{verlinde}
\def\theenumi{\roman{enumi}}
\def\labelenumi{(\theenumi)}
\begin{enumerate}
\item
We have
$$
\dim H^0(\SVB(C_{\eta}),\Theta_{\SVB}(C_{\eta})^{\kappa})=
(n(n+\kappa)^{n-1})^{g-1}\cdot
\sum_{\mu\in \rho+P_{\kappa}}|J(\rho,\mu)|^{2(1-g)}
$$
\item
For $(a',b')\in\SA'$ we have
$$
\dim H^0(\SPB,\Theta_{\SPB}^{\kappa}(a',b'))=
(n(n+\kappa)^{n-1})^{g-2}\cdot
\sum_{\mu\in \rho+P_{\kappa}}|J(\lambda,\mu)|^2|J(\rho,\mu)|^{2(1-g)}
\quad,
$$
where $\lambda=\rho+\sum_{i=1}^na'_{n-i+1}\epsilon_i$.
\end{enumerate}
\end{proposition}

\begin{proof}
This is a well established fact, only the shape of
the formulas is maybe a bit unusual.
To explain the formulas we will employ the notation of \cite{B}.
Let $X$ be a smooth projective curve over $\C$, let
$x_1,\dots,x_m$ be distinct points on $X$ and let 
$\lambda_1,\dots,\lambda_m$ be dominant weights of the
Lie-algebra $\sln$. To these data there is associated a finite dimensional
vector space $V_{X}((x_1,\dots,x_m),(\lambda_1,\dots,\lambda_m))$,
called the space of conformal blocks.
As shown in \cite{BL} and \cite{Pauly}, we have
\begin{eqnarray*}
\dim H^0(\SVB(C_{\eta}),\Theta_{\SVB}(C_{\eta})^{\kappa}) &=&
\dim V_{Y}(\emptyset)
\\
\dim H^0(\SPB,\Theta_{\SPB}^{\kappa}(a',b')) &=&
\dim V_{\Ct}((p_1,p_2),(\lambda,\lambda^*))
\end{eqnarray*}
where $Y$ is some smooth projective curve of genus $g$,
$\lambda$ is the dominant weight
$
\sum_{i=1}^na'_{n-i+1}\epsilon_i
$
and
$
\lambda^*=-w_0\lambda=\sum_{i=1}^nb'_{n-i+1}\epsilon_i
$
is its dual.

By \cite{B}, Cor. 9.8 we have
\begin{eqnarray*}
\dim V_{Y}(\emptyset) &=&
(n(n+\kappa)^{n-1})^{g-1}\cdot
 \sum_{\mu\in P_{\kappa}}\frac{1}{\Delta(t_{\mu})^{g-1}}
\\
\dim V_{\Ct}((p_1,p_2),(\lambda,\lambda^*))&=&
(n(n+\kappa)^{n-1})^{g-2}\cdot
 \sum_{\mu\in P_{\kappa}}
\frac
{\Tr_{V_{(\lambda,\lambda^*)}}(t_{\mu})}
{\Delta(t_{\mu})^{g-2}}
\end{eqnarray*}
where 
$
\Delta(t_{\mu})=|J(\rho,\rho+\mu)|^2
$
and
$
\Tr_{V_{(\lambda,\lambda^*)}}(t_{\mu})=
|J(\rho+\lambda,\rho+\mu)|^2/\Delta(t_{\mu})
$.
\end{proof}

From Proposition \ref{verlinde} it is immediate that
for the proof of
equality $(2)$ it is sufficient to show the following
lemma, which is elementary in its statement 
but which I could not prove without the help of Don Zagier.
\begin{lemma}
\label{zagier}
For every
$\mu\in \rho+P_{\kappa}$ we have
$$
\sum_{\lambda\in\rho+P_{\kappa}}
\gamma(\lambda)\cdot
|J(\lambda,\mu)|^2=
n(n+\kappa)^{n-1}
\quad,
\eqno(3)
$$
where for $\lambda=\sum\lambda_i\epsilon_i\in P$ we set
$\gamma(\lambda):=\kappa+n-1-\lambda_1$.
\end{lemma}

\begin{proof}
(The proof of this lemma is due to Don Zagier).
Let $m:=\kappa+n$ and let $\zeta_m:=\exp(2\pi i/m)$.
Observe that the mapping 
$
\lambda=
\sum_i\lambda_i\epsilon_i
\mapsto
\{\lambda_1,\dots,\lambda_n\}
$
is a bijection from the set $\rho+P_{\kappa}$ to
the set $Q_0$ of all subsets $A$ of $N:=\{0,\dots,m-1\}$
with $|A|=n$ and $0\in A$.
It follows directly from the definitions that
if $\lambda,\mu\in\rho+P_{\kappa}$ are mapped
to $A,B\in Q_0$, then we have
$$
|J(\lambda,\mu)|^2=|\Delta_{A,B}|^2
\quad,
$$
where $\Delta_{A,B}$ is the determinant of the
$n\times n$-sub-matrix of the matrix 
$
M_m=(\zeta_m^{a,b})_{0\leq a,b\leq m-1}
$
corresponding to the rows and columns with indices in the
sets $A$ and $B$.
Thus the assertion of the lemma holds if and only if
$$
\sum_{A\in Q_0}
\gamma(A)|\Delta_{A,B}|^2
=
(m-n)m^{n-1}
\qquad\qquad
(B\in Q)
\eqno(4)
$$
where $\gamma(A)=n-\max(A)-1$ and $Q$ is the set of all subsets
$B\subset N$ with $|B|=n$. (The condition ``$0\in B$'' can be omitted
since a translation of the set $B$ just multiplies every determinant
$\Delta_{A,B}$ by a root of unity and does not have any influence
on the expression $|\Delta_{A,B}|^2$.)
We can further optically simplify the formula by replacing $M_m$ by
$M^*_m=m^{-1/2}M_m$ and thus $\Delta_{AB}$ by 
$\Delta^*_{AB}=m^{-n/2}\Delta_{AB}$ 
(which is reasonable since $M_m\overline{M_m^t}=m\cdot\Id$
and therefore $M^*_m$ is unitary):
$$
\sum_{A\in Q_0}
\gamma(A)|\Delta^*_{AB}|^2 
=
\frac{m-n}{m}
\qquad\qquad
(B\in Q).
\eqno(5)
$$

We compute first the left hand side of (5)
leaving out the factor $\gamma(A)$ and the condition ``$0\in A$''.
The numbers $\{\Delta^*_{AB}\}_{A,B\in Q}$ are nothing else but
the matrix coefficients of the $n$-th exterior product
$\bigwedge^n(M^*_m)$ of the operator represented by the matrix
$M^*_m$.
The fact that $M^*_m$ is unitary remains true also for $\bigwedge^n$
of this matrix; 
consequently we have 
$
\bigwedge^n(M^*_m)\bigwedge^n(\overline{M^{*t}_m})=\Id_{\binom{m}{n}}
$
or explicitly
$
\sum_{A\in Q}\Delta^*_{AB}\overline{\Delta^*_{AB'}}=\delta_{BB'}
$
for $B,B'\in Q$.
The special case $B=B'$ of this yields
$$
\sum_{A\in Q}|\Delta^*_{AB}|^2 = 1,
\qquad\qquad
(B\in Q).
\eqno(6)
$$

It remains to show that the left hand side of the equation (5)
differs by the factor $(m-n)/m$ from the left hand side of the 
equation (6).
For this we give yet another description of the set $Q$.
We denote by $\Qb$ the set of all mappings
$\alpha:\Z/n\to\Z/m$ that are 
``cyclically strictly decreasing''
(i.e. they can be lifted to a mapping $a:\Z\to\Z$ 
for which $a(i)>a(i+1)>\dots>a(i+n)=a(i)-m$ for all $i\in\Z$ holds).
Let $\Qb_0\subset\Qb$ be the subset defined by the additional
property $\alpha(0)=0$.
We have a diagram of mappings between sets as follows:
$$
\xymatrix@R=2ex{
\text{$\Qb_0$} \ar[r] \ar@{^(->}[d] &
Q_0 \ar@{^(->}[d]
\\
\text{$\Qb$} \ar[r] &
Q
}
$$
where the map
$\Qb\to Q$
maps $\alpha$ to the subset $A$ of $N$ which 
corresponds to the image of $\alpha$  via the obvious bijection
$N\isomto \Z/m$
and $\Qb_0\to Q_0$ is the restriction of
$\Qb\to Q$
to the subset $\Qb_0$.
There is an operation of the cyclic group $C_n=\Z/n$ on
$\Qb$ defined by $\alpha\mapsto\alpha(\cdot+j)$ and 
the mapping $\Qb\to Q$ identifies $Q$ with the quotient 
$\Qb/C_n$.
On the other hand, there is also an operation of the
cyclic group $C_m=\Z/m$ on $\Qb$, given by $\alpha\mapsto\alpha +k$
and $\Qb_0$ is a system of representatives for this operation.
The mapping $\Qb_0\to Q_0$ is a bijection.

For $\alpha\in\Qb$ with lifting $a:\Z\to\Z$ and $i\in \Z$
the number 
$\gamma_i(\alpha):=a(i)-a(i+1)-1$
depends only on $\alpha$ and on $i$ mod $n$;
clearly if $\alpha_0$ is in $\Qb$ and $A_0$ the
corresponding element in $Q$, then 
$
\sum_{\alpha\mapsto A_0}\gamma_0(\alpha)
=
\sum_{\text{$i$ mod $n$}}\gamma_i(\alpha_0)=m-n
$.
We write $\Delta^*_{\alpha B}:=\Delta^*_{AB}$ where $A\in Q$
is the element corresponding to $\alpha$.
Equation (5) can now be seen as follows:
\begin{eqnarray*}
\sum_{A\in Q_0}
\gamma(A)|\Delta^*_{AB}|^2
&=&
\sum_{\alpha\in\Qb_0}
\gamma_0(\alpha)|\Delta^*_{\alpha B}|^2
=
\\
&\overset{(*)}{=}&
\frac{1}{m}
\sum_{\alpha\in\Qb}
\gamma_0(\alpha)|\Delta^*_{\alpha B}|^2 
=
\frac{1}{m}
\sum_{A\in Q}
|\Delta^*_{AB}|^2
\sum_{\alpha\mapsto A}
\gamma_0(\alpha)
=
\\
&=&
\frac{1}{m}
\sum_{A\in Q}
(m-n)|\Delta^*_{AB}|^2
\\
&\overset{(6)}{=}&
\frac{m-n}{m}
\quad.
\end{eqnarray*}
Here $(*)$ follows from the fact that 
$
\gamma_0(\alpha+k)|\Delta^*_{\alpha+k,B}|=
\gamma_0(\alpha)|\Delta^*_{\alpha,B}|
$
for $k\in C_m$.

%weiter
\end{proof}

\end{document}